\documentclass{article}
\usepackage{psfig}
\usepackage{graphicx}
\usepackage{lscape}

\def\qed{\hfill {\hbox{${\vcenter{\vbox{               
   \hrule height 0.4pt\hbox{\vrule width 0.4pt height 6pt
   \kern5pt\vrule width 0.4pt}\hrule height 0.4pt}}}$}}}

\newtheorem{theorem}{Theorem}
\newtheorem{definition}{Definition}

\newtheorem{conjecture}{Conjecture}

{\qed\par\medskip}

\title{\Large \textbf{MEANDER KNOTS AND LINKS}}

\author{
Slavik Jablan\footnote{The Mathematical Institute, Belgrade, 11000, Serbia, sjablan@gmail.com}\\
Ljiljana Radovi\' c\footnote{The Faculty of Mechanical Engineering, Ni\v s, 18000, Serbia, ljradovic@gmail.com}}

\begin{document}

\maketitle

\begin{abstract}

We introduced concept of meander knots, 2-component meander links and multi-component meander links and derived different families of meander knots and links from open meanders with $n\le 16$ crossings. We also defined semi-meander knots (or knots with ordered Gauss code) and their product.
\end{abstract}


\section{Introduction}

Meanders are combinatorial objects with some topological properties of the interplay between planarity and connectedness. They
correspond to the systems formed by the intersections of two curves in the
plane, with equivalence up to homeomorphism within the plane. They occur in polymer physics, algebraic geometry, mathematical theory of mazes, and the study of
planar algebras, especially the Temperley-Lieb algebra. As the main source to the theory of meanders we used the paper \cite{croix}. For applications of
meanders, the reader is referred to \cite{1,2,3}.

\begin{definition}
An open meander is a configuration consisting of an oriented simple curve, and a line in the plane, the axis of the meander,  that cross a finite number of times and intersect only transversally. Two open meanders are equivalent if there is a homeomorphism of the plane that maps one meander to the other.
\end{definition}

The  {\it order of a meander} is the number of crossings between the meander curve and meander axis. Since a line and a simple curve are homeomorphic, and we are only interested in defining meanders up to homeomorphism, their roles can be reversed by imposing an orientation on the line instead of the curve. However, in the enumeration of meanders we will always distinguish the meander curve from the meander line -- the axis. Usually, meanders are classified according to their order. One of the main problems the theory of meanders is their enumeration.

A open meander curve and meander axis have two loose ends each. Depending on the number of crossings, the loose ends of the meander curve belong to different half-planes defined by the axis for open meanders with an odd order, and to same half-plane for meanders with an even number of crossings. We will first consider the first case. In this case we are able to make a closure of the meander: to join each of loose ends of the meander curve with an end of the axis, and obtain a knot or link shadow that will be called {\it meandric knot shadow}. Because for making a closure we have two possibilities, we will always choose this one producing meandric knot shadow without loops. After that, by introducing undercrossings and overcrossings along the meander knot shadow axis, we will turn it into a knot diagram.

\begin{definition}
An alternating knot which has a minimal meander diagram is called \it{meander knot}.
\end{definition}

For representing meanders we will use their corresponding arch configurations.

\begin{definition}
An arch configuration is a planar configuration consisting of pairwise non-intersecting semicircular (or piecewise-linear) arches lying on the same side of an oriented line, arranged such that the feet of the arches are equally spaced along the line.
\end{definition}

Arch configurations play an essential role in the enumerative theory of meanders. They are used to obtain the canonical representatives of meanders that are used in subsequent constructions. They also have a natural link to the Temperley-Lieb algebra. A meandric system is the superposition of an ordered pair
of arch configurations of the same order, with the first configuration as the upper configuration, and the second as the lower configuration.

\begin{figure}[th]
\centerline{\psfig{file=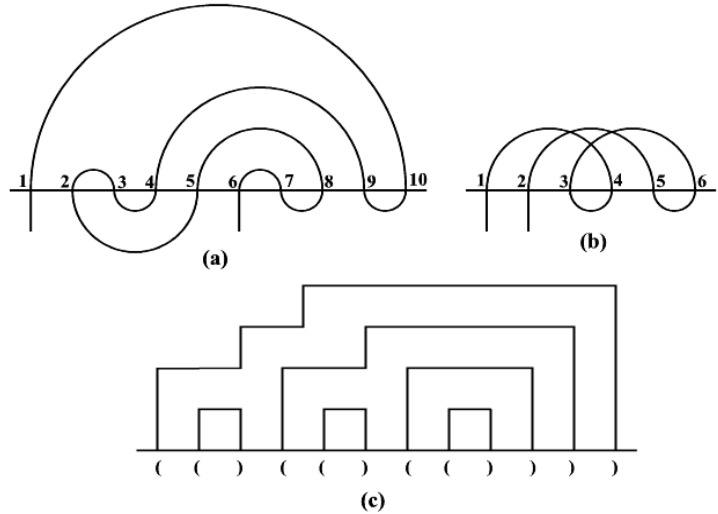,width=3.8in}} \vspace*{8pt}
\caption{(a) Open meander given by meander permutation $(1,10,9,4,3,2,5,8,7,6)$; (b) non-realizable sequence $(1,4,3,6,5,2)$; (c) piecewise-linear upper arch configuration given by Dyck word $(()(()(())))$.
\label{m1}}
\end{figure}

The problem of enumerating meanders has a long history, with interest dating back at least to work by Poincar\' e on differential geometry, though
the modern study appears to have been inspired by Arnol'd in \cite{1}. If the intersections along the axis are enumerated by $1,2,3,...,n$ every open meander can be described by a meander permutation of order $n$: the sequence of $n$ numbers describing the path of the meander curve. For example, the open meander (Fig. \ref{m1}a) is coded by the meander permutation $(1,10,9,4,3,2,5,8,7,6).$ Hence, enumeration of open meanders is based on the derivation of meander permutations. Meander permutations play an important role in the mathematical theory of mazes \cite{4}. In every meander permutation alternate odd and even numbers, i.e., to the upper and lower arch configuration correspond sequences of opposite parity. However, this condition does not completely characterize meander permutations. E.g., the sequence $(1,4,3,6,5,2)$ is not a meander permutation, since it not satisfies planarity condition: arc $(1,4)$ and $(3,6)$ will produce a nugatory crossing (Fig. \ref{m1}b). Hence, the most important property of meander permutations is that all arcs must be nested, in order to not produce nugatory crossings. Among different techniques to achieve this, the fastest algorithms for deriving meanders are based on encoding of arch configurations as words in the Dyck language: the language of balanced parentheses on the alphabet $\{x,y\}$ which can be generated by the production rules:

$$s\rightarrow e, \quad s \rightarrow xsys$$

\noindent where $e$ is empty word. Arch configurations of order $n$ are in bijective correspondence with Dyck words of length $2n$. E.g., to arch configuration from Fig. \ref{m1}c corresponds the Dyck word $(()(()(())))$. Every open meander can be uniquely described by Dyck words corresponding to upper and lower arch configuration, with two additional symbols 1 for positions of loose ends. E.g., the open meander coded by the meander permutation $(1,10,9,4,3,2,5,8,7,6)$ (Fig. \ref{m1}a) is defined by upper and lower Dyck word: $\{(()((()))),1(())1()()\}$. Algorithms for meander enumeration can be found in different sources (e.g., \cite{en1,en2}). For the derivation of open meanders we used the {\it Mathematica} program "Open meanders" by David Bevan ({\tt http://demonstrations.wolfram.com/OpenMeanders/}) \cite{bev}, which we modified in order to compute open meander permutations instead of Dyck words.

Obtained meander knots and links will be given by their Gauss codes, Dowker-Thisllethwaite codes, and Conway symbols \cite{con,rol,jab}. All computations are made in the program "LinKnot" \cite{jab}.

\section{Meander knots}

\begin{figure}[th]
\centerline{\psfig{file=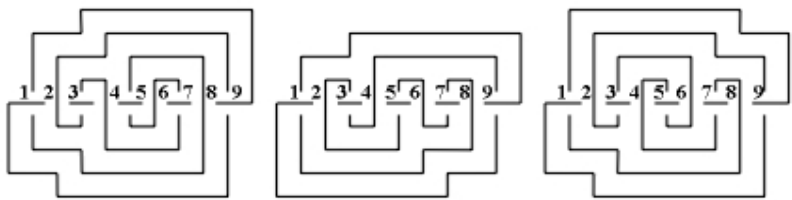,width=4.2in}} \vspace*{8pt}
\caption{Non-isomorphic minimal meander diagrams of the knot $9_7=3\,4\,2$.
\label{m2}}
\end{figure}

Gauss codes of alternating meander knot diagrams can be obtained if to the sequence $1,2,\ldots n$ we add a meander permutation of order $n$ (and/or its reverse), where $n$ is an odd number and in the obtained sequence alternate the signs of successive numbers. E.g., from meander permutation $(1,8,5,6,7,4,3,2,9)$ we obtain Gauss code $$\{-1,2,-3,4,-5,6,-7,8,-9,1,-8,5,-6,7,-4,3,-2,9\}$$
\noindent which corresponds to rational alternating knot $9_7$ given by the Conway symbol $3\,4\,2$. The same knot can be obtained from meander permutations $(1,8,7,6,5,2,$ $3,4,9)$ and $(1,8,7,4,5,6,3,2,9)$, giving Gauss codes
$$\{-1,2,-3,4,-5,6,-7,8,-9,1,-8,7,-6,5,-2,3,-4,9\}$$
\noindent and
$$\{-1,2,-3,4,-5,6,-7,8,-9,1,-8,7,-4,5,-6,3,-2,9\},$$
\noindent which can be converted to Dowker-Thistlethwaite codes $\{\{9\},\{4,12,16,18,14,$ $2,10,8,6\}\}$,  $\{\{9\},\{4,12,16,18,14,2,10,8,6\}\}$, and $\{\{9\},\{4,12,18,16,14,2,10,$ $6,8\}\}$, corresponding to 3 non-isomorphic minimal diagrams of the knot $9_7$ $=3\,4\,2$ (Fig. \ref{m2}). The reverse of a meander permutation we add only if the absolute value of the last number of the first sequence is different from the absolute value of the first number of the second sequence, in order to avoid curls. However, if an alternating knot has an alternating minimal meander diagram, all its minimal diagrams are not necessarily meander diagrams. E.g., knot $K11a192=3\,1\,2\,1\,2\,2$ has 10 non-isomorphic flype-equivalent minimal diagrams, but only two of them are meander diagrams. Hence, an alternating knot will be called {\it meander knot} if it has at least one minimal meander diagram. Every alternating meander diagram is a positive knot diagram. From every alternating diagram, by crossing changes we obtain (minimal and non-minimal) meander diagrams of different knots. Among meander knots we can also obtain non-prime meander knots, so our discussion we will restrict only to prime meander knot diagrams.

\begin{figure}[th]
\centerline{\psfig{file=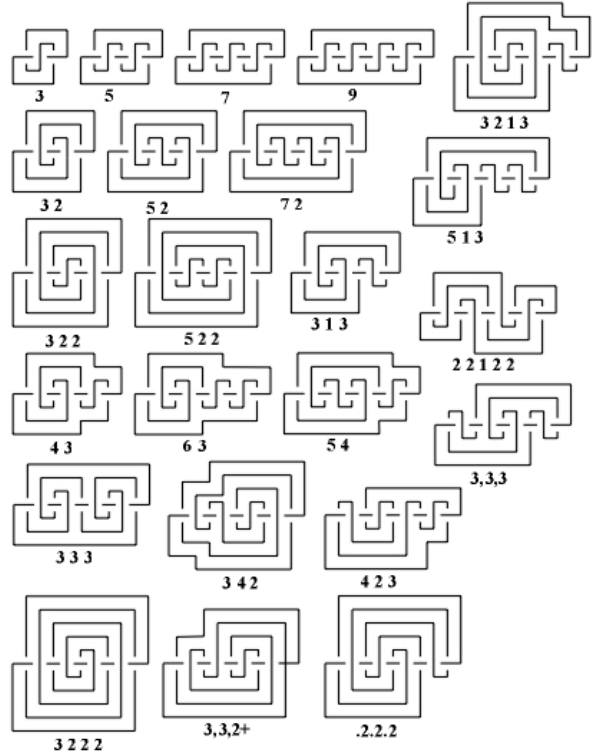,width=3.2in}} \vspace*{8pt}
\caption{Alternating meander knots with at most $n=9$ crossings.
\label{alt}}
\end{figure}

\begin{figure}[th]
\centerline{\psfig{file=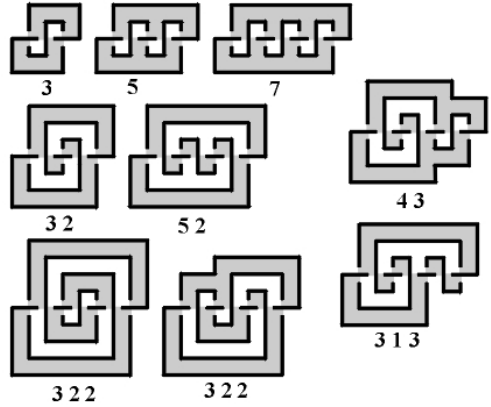,width=2.6in}} \vspace*{8pt}
\caption{Checker-board colorings of meander diagrams of meander knots with at most $n=7$ crossings.
\label{bw}}
\end{figure}

The natural question which arises is: find all alternating meander knots with $n$ crossings, where $n$ is an odd number. In the following table is given the number of open meanders $OM$ with $n$ crossings and the number $AMK$ of their corresponding alternating meander knots.

\bigskip

\begin{tabular}{|c|c|c|}   \hline
$n$ & $OM$ & $AMK$  \\  \hline
1 & 1 &  1  \\  \hline
3 & 2 &  1  \\  \hline
5 & 8 &   2 \\  \hline
7 & 42 &   5  \\  \hline
9 & 262 &   15  \\  \hline
11 & 1828 &  52  \\  \hline
13& 13820 &  233   \\  \hline
15& 110954 &  1272   \\  \hline
\end{tabular}

\bigskip

Alternating meander knots with at most $n=9$ crossings are illustrated in Fig. \ref{alt}. Notice that alternating meander knots are organized in families. E.g., the family $3_1=3$, $5_1=5$, $7_1=7$, $9_1=9$, $K11a367=11$, $\ldots $, or the family of twist knots $5_2=3\,2$, $7_2=5\,2$, $9_2=7\,2$, $K11a247=9\,2$, $\ldots $. Since all arcs have an integer radius, for every alternating meander knot we can define its height: minimal vertical dimension of its minimal meander diagrams. This dimension is the same for all meander knots belonging to a family. Every meander diagram produces very nice checkerboard coloring (Fig. \ref{bw}).

In the next table are given all alternating meander prime knots with at most $n=11$ crossing. Every knot is given by its standard symbol from Rolfsen book \cite{rol} or ordering number from "Knotscape", Conway symbol, Dowker-Thistlethwaite code of its meander diagram and short Gauss code. Short Gauss code means that Gauss code of every meander knot is beginning with the sequence $(-1)^ii$ ($1\le i\le n$), followed by meander permutation part, so the first sequence is omitted.

\bigskip

\tiny

\noindent \begin{tabular}{|c|c|c|c|}   \hline
$3_1$ & $3$ & $\{\{3\},\{4,6,2\}\}$ &  $\{1,-2,3\}$ \\ \hline
$5_1$ & 5 & $\{\{5\},\{6,8,10,2,4\}\}$ &  $\{1,-2,3,-4,5\}$ \\ \hline
$5_2$ & $3\,2$ & $\{\{5\},\{4,8,10,2,6\}\}$ & $\{3,-2,1,-4,5\}$ \\ \hline
$7_1$ & 7 & $\{\{7\},\{8,10,12,14,2,4,6\}\}$ & $\{1,-2,3,-4,5,-6,7\}$ \\ \hline
$7_2$ & $5\,2$ & $\{\{7\},\{4,10,14,12,2,8,6\}\}$ & $\{1,-2,7,-6,5,-4,3\}$ \\ \hline
$7_3$ & $4\,3$ & $\{\{7\},\{6,10,12,14,2,4,8\}\}$ & $\{3,-2,1,-4,5,-6,7\}$ \\ \hline
$7_4$ & $3\,1\,3$ & $\{\{7\},\{6,10,12,14,4,2,8\}\}$ & $\{1,-4,3,-2,7,-6,5\}$ \\ \hline
$7_5$ & $3\,2\,2$ & $\{\{7\},\{4,10,14,12,2,6,8\}\}$ & $\{3,-4,5,-2,1,-6,7\}$ \\ \hline
$9_1$ & 9 & $\{\{9\},\{10,12,14,16,18,2,4,6,8\}\}$ & $\{1,-2,3,-4,5,-6,7,-8,9\}$ \\ \hline
$9_2$ & $7\,2$ & $\{\{9\},\{4,12,18,16,14,2,10,8,6\}\}$ & $\{1,-2,9,-8,7,-6,5,-4,3\}$ \\ \hline
$9_3$ & $6\,3$ & $\{\{9\},\{8,12,14,16,18,2,4,6,10\}\}$ & $\{3,-2,1,-4,5,-6,7,-8,9\}$ \\ \hline
$9_4$ & $5\,4$ & $\{\{9\},\{6,12,14,18,16,2,4,10,8\}\}$ & $\{5,-4,3,-2,1,-6,7,-8,9\}$  \\ \hline
$9_5$ & $5\,1\,3$ & $\{\{9\},\{6,12,14,18,16,4,2,10,8\}\}$ & $\{1,-4,3,-2,9,-8,7,-6,5\}$ \\ \hline
$9_6$ & $5\,2\,2$ & $\{\{9\},\{4,12,18,14,16,2,6,8,10\}\}$ & $\{3,-4,5,-6,7,-2,1,-8,9\}$ \\ \hline
$9_7$ & $3\,4\,2$ & $\{\{9\},\{4,12,18,16,14,2,8,10,6\}\}$ & $\{1,-2,9,-8,7,-6,3,-4,5\}$  \\ \hline
$9_9$ & $4\,2\,3$ & $\{\{9\},\{6,12,14,18,16,2,4,8,10\}\}$ & $\{3,-4,5,-2,1,-6,7,-8,9\}$  \\ \hline
$9_{10}$ & $3\,3\,3$ & $\{\{9\},\{8,12,14,16,18,4,2,6,10\}\}$ & $\{1,-2,3,-6,5,-4,9,-8,7\}$  \\ \hline
$9_{13}$ & $3\,2\,1\,3$ & $\{\{9\},\{6,12,14,16,18,4,2,10,8\}\}$ & $\{1,-6,3,-4,5,-2,9,-8,7\}$ \\ \hline
$9_{16}$ & $3,3,2+$ & $\{\{9\},\{4,12,16,18,14,2,8,10,6\}\}$ & $\{5,-6,7,-2,3,-4,1,-8,9\}$ \\ \hline
$9_{18}$ & $3\,2\,2\,2$ & $\{\{9\},\{4,12,18,14,16,2,6,10,8\}\}$ & $\{5,-4,3,-6,7,-2,1,-8,9\}$  \\ \hline
$9_{23}$ & $2\,2\,1\,2\,2$ & $\{\{9\},\{4,10,18,14,2,16,6,12,8\}\}$ & $\{7,-6,1,-2,5,-8,9,-4,3\}$ \\ \hline
$9_{35}$ & $3,3,3$ & $\{\{9\},\{8,12,16,14,18,4,2,6,10\}\}$ & $\{3,-2,1,-6,5,-4,9,-8,7\}$ \\ \hline
$9_{38}$ & $.2.2.2$ & $\{\{9\},\{6,10,14,18,4,16,2,8,12\}\}$ & $\{1,-6,5,-2,9,-8,3,-4,7\}$ \\ \hline
\end{tabular}

\noindent \begin{tabular}{|c|c|c|c|}   \hline
$K11a94$ & $2\,2\,1,3,2+$ & $\{\{11\},\{4,10,22,16,2,20,18,6,12,14,8\}\}$ & $\{7,-8,9,-6,1,-2,5,-10,11,-4,3\}$  \\ \hline
$K11a95$ & $2\,4\,1\,2\,2$ & $\{\{11\},\{4,10,22,16,2,20,18,6,14,12,8\}\}$ & $\{9,-8,7,-6,1,-2,5,-10,11,-4,3\}$ \\ \hline
$K11a186$&$2\,2\,3\,2\,2$&$\{\{11\},\{4,12,14,16,20,2,10,6,22,8,18\}\}$&$\{9,-8,1,-2,5,-6,7,-10,11,-4,3\}$\\ \hline
$K11a186$&$2\,2\,3\,2\,2$&$\{\{11\},\{4,12,14,16,20,2,10,6,22,8,18\}\}$&$\{9,-8,1,-2,5,-6,7,-10,11,-4,3\}$\\ \hline
$K11a191$&$4\,2\,1\,2\,2$&$\{\{11\},\{4,12,14,18,20,2,10,22,6,8,16\}\}$&$\{7,-6,1,-2,5,-8,9,-10,11,-4,3\}$\\ \hline
$K11a192$&$3\,1\,2\,1\,2\,2$&$\{\{11\},\{4,12,14,18,20,2,10,22,8,6,16\}\}$&$\{9,-8,1,-4,3,-2,7,-10,11,-6,5\}$\\ \hline
$K11a200$&$2\,2,2\,2,3$&$\{\{11\},\{4,12,16,14,20,2,10,6,22,8,18\}\}$&$\{9,-8,1,-2,7,-6,5,-10,11,-4,3\}$\\ \hline
$K11a234$&$7\,2\,2$&$\{\{11\},\{4,14,16,18,20,22,2,12,6,8,10\}\}$&$\{3,-4,5,-6,7,-8,9,-2,1,-10,11\}$\\ \hline
$K11a235$&$3\,4\,2\,2$&$\{\{11\},\{4,14,16,18,20,22,2,12,6,10,8\}\}$&$\{5,-4,3,-6,7,-8,9,-2,1,-10,11\}$\\ \hline
$K11a236$&$3\,2\,2\,2\,2$&$\{\{11\},\{4,14,16,18,20,22,2,12,10,8,6\}\}$&$\{5,-6,7,-4,3,-8,9,-2,1,-10,11\}$\\ \hline
$K11a237$&$2\,2,3,3+$&$\{\{11\},\{4,14,16,18,22,20,2,12,8,6,10\}\}$&$\{1,-2,11,-10,3,-6,5,-4,9,-8,7\}$\\ \hline
$K11a238$&$5\,2\,2\,2$&$\{\{11\},\{4,14,16,20,18,22,2,12,10,8,6\}\}$&$\{1,-2,11,-10,3,-4,9,-8,7,-6,5\}$\\ \hline
$K11a240$&$5,3,2+$&$\{\{11\},\{4,14,18,20,22,16,2,10,12,6,8\}\}$&$\{5,-6,7,-8,9,-2,3,-4,1,-10,11\}$\\ \hline
$K11a241$&$3\,2,3,2+$&$\{\{11\},\{4,14,18,20,22,16,2,10,12,8,6\}\}$&$\{7,-6,5,-8,9,-2,3,-4,1,-10,11\}$\\ \hline
$K11a242$&$5\,4\,2$&$\{\{11\},\{4,14,18,20,22,16,2,12,10,6,8\}\}$&$\{5,-6,7,-8,9,-4,3,-2,1,-10,11\}$\\ \hline
$K11a243$&$3\,2\,4\,2$&$\{\{11\},\{4,14,18,20,22,16,2,12,10,8,6\}\}$&$\{1,-2,11,-10,9,-8,3,-4,7,-6,5\}$\\ \hline
$K11a244$&$6^*2\,1.2:2\,0.2\,0$&$\{\{11\},\{4,14,18,22,16,20,2,8,12,6,10\}\}$&$\{7,-4,3,-8,9,-2,5,-6,1,-10,11\}$\\ \hline
$K11a245$&$3,3,2+++$&$\{\{11\},\{4,14,20,22,16,18,2,12,10,8,6\}\}$&$\{1,-10,9,-8,5,-6,7,-2,3,-4,11\}$\\ \hline
$K11a246$&$3\,6\,2$&$\{\{11\},\{4,14,20,22,18,16,2,12,10,8,6\}\}$&$\{1,-2,11,-10,9,-8,7,-6,3,-4,5\}$\\ \hline
$K11a247$&$9\,2$&$\{\{11\},\{4,14,22,20,18,16,2,12,10,8,6\}\}$&$\{1,-2,11,-10,9,-8,7,-6,5,-4,3\}$\\ \hline
$K11a263$&$3,3,3,2$&$\{\{11\},\{6,8,16,2,20,22,18,4,12,14,10\}\}$&$\{7,-8,9,-4,5,-6,1,-2,3,-10,11\}$\\ \hline
$K11a291$&$6^*2.4:2\,0$&$\{\{11\},\{6,10,16,22,4,18,20,2,8,12,14\}\}$&$\{1,-8,7,-2,11,-10,3,-4,5,-6,9\}$\\ \hline
$K11a292$&$6^*2.3\,1:2\,0$&$\{\{11\},\{6,10,16,22,4,18,20,2,8,14,12\}\}$&$\{1,-8,7,-2,11,-10,3,-6,5,-4,9\}$\\ \hline
$K11a298$&$6^*3\,1.2:2\,0$&$\{\{11\},\{6,10,18,2,16,20,22,4,8,14,12\}\}$&$\{1,-8,5,-6,7,-2,11,-10,3,-4,9\}$\\ \hline
$K11a299$&$6^*4.2:2\,0$&$\{\{11\},\{6,10,18,2,16,22,20,4,8,14,12\}\}$&$\{1,-6,5,-2,11,-10,9,-8,3,-4,7\}$\\ \hline
$K11a318$&$6^*2.2.2.2\,0.2\,0$&$\{\{11\},\{6,12,16,22,18,2,20,4,8,14,10\}\}$&$\{7,-6,1,-2,5,-8,9,-4,3,-10,11\}$\\ \hline
$K11a319$&$6^*3.2.2:2\,0$&$\{\{11\},\{6,12,16,22,18,4,20,2,8,10,14\}\}$&$\{1,-8,7,-2,3,-4,11,-10,5,-6,9\}$\\ \hline
$K11a320$&$6^*2.2.3.2\,0$&$\{\{11\},\{6,12,16,22,18,4,20,2,10,8,14\}\}$&$\{1,-6,5,-2,11,-10,3,-4,9,-8,7\}$\\ \hline
$K11a334$&$5\,2\,4$&$\{\{11\},\{6,14,16,18,20,22,2,4,12,8,10\}\}$&$\{3,-4,5,-6,7,-2,1,-8,9,-10,11\}$\\ \hline
$K11a335$&$4\,2\,2\,3$&$\{\{11\},\{6,14,16,18,20,22,2,4,12,10,8\}\}$&$\{5,-4,3,-6,7,-2,1,-8,9,-10,11\}$\\ \hline
$K11a336$&$5\,2\,1\,3$&$\{\{11\},\{6,14,16,18,20,22,4,2,12,8,10\}\}$&$\{1,-8,3,-4,5,-6,7,-2,11,-10,9\}$\\ \hline
$K11a337$&$3\,2\,2\,1\,3$&$\{\{11\},\{6,14,16,18,20,22,4,2,12,10,8\}\}$&$\{1,-8,3,-4,7,-6,5,-2,11,-10,9\}$\\ \hline
$K11a338$&$4,3,3+$&$\{\{11\},\{6,14,16,20,22,18,2,4,10,12,8\}\}$&$\{5,-6,7,-2,3,-4,1,-8,9,-10,11\}$\\ \hline
$K11a339$&$4\,4\,3$&$\{\{11\},\{6,14,16,20,22,18,2,4,12,10,8\}\}$&$\{5,-6,7,-4,3,-2,1,-8,9,-10,11\}$\\ \hline
$K11a340$&$3\,1,3,3+$&$\{\{11\},\{6,14,16,20,22,18,4,2,10,12,8\}\}$&$\{1,-8,5,-6,7,-2,3,-4,11,-10,9\}$\\ \hline
$K11a341$&$3\,4\,1\,3$&$\{\{11\},\{6,14,16,20,22,18,4,2,12,10,8\}\}$&$\{1,-8,7,-6,3,-4,5,-2,11,-10,9\}$\\ \hline
$K11a342$&$7\,4$&$\{\{11\},\{6,14,16,22,20,18,2,4,12,10,8\}\}$&$\{1,-2,3,-4,11,-10,9,-8,7,-6,5\}$\\ \hline
$K11a343$&$7\,1\,3$&$\{\{11\},\{6,14,16,22,20,18,4,2,12,10,8\}\}$&$\{1,-4,3,-2,11,-10,9,-8,7,-6,5\}$\\ \hline
$K11a353$&$6^*3.2:2\,0.2\,0$&$\{\{11\},\{8,12,16,18,22,4,20,2,6,10,14\}\}$&$\{1,-4,7,-8,3,-2,9,-6,5,-10,11\}$\\ \hline
$K11a354$&$6^*2.2:3.2 0$&$\{\{11\},\{8,12,16,20,2,18,22,6,4,10,14\}\}$&$\{3,-2,1,-8,7,-4,11,-10,5,-6,9\}$\\ \hline
$K11a355$&$6\,2\,3$&$\{\{11\},\{8,14,16,18,20,22,2,4,6,12,10\}\}$&$\{3,-4,5,-2,1,-6,7,-8,9,-10,11\}$\\ \hline
$K11a356$&$3\,3\,2\,3$&$\{\{11\},\{8,14,16,18,20,22,2,6,4,12,10\}\}$&$\{1,-2,3,-8,5,-6,7,-4,11,-10,9\}$\\ \hline
$K11a357$&$3\,2\,1\,2\,3$&$\{\{11\},\{8,14,16,18,20,22,6,4,2,12,10\}\}$&$\{1,-6,3,-4,5,-2,11,-10,7,-8,9\}$\\ \hline
$K11a358$&$6\,5$&$\{\{11\},\{8,14,16,18,22,20,2,4,6,12,10\}\}$&$\{5,-4,3,-2,1,-6,7,-8,9,-10,11\}$\\ \hline
$K11a359$&$5\,3\,3$&$\{\{11\},\{8,14,16,18,22,20,2,6,4,12,10\}\}$&$\{1,-2,3,-6,5,-4,11,-10,9,-8,7\}$\\ \hline
$K11a360$&$5\,1\,2\,3$&$\{\{11\},\{8,14,16,18,22,20,6,4,2,12,10\}\}$&$\{1,-6,3,-4,5,-2,11,-10,9,-8,7\}$\\ \hline
$K11a361$&$3\,2,3,3$&$\{\{11\},\{8,14,18,16,20,22,4,2,6,12,10\}\}$&$\{5,-4,1,-2,3,-8,7,-6,11,-10,9\}$\\ \hline
$K11a362$&$5,3,3$&$\{\{11\},\{8,14,18,16,22,20,4,2,6,12,10\}\}$&$\{3,-2,1,-6,5,-4,11,-10,9,-8,7\}$\\ \hline
$K11a363$&$5\,1\,5$&$\{\{11\},\{8,16,14,18,22,20,6,4,2,12,10\}\}$&$\{1,-6,5,-4,3,-2,11,-10,9,-8,7\}$\\ \hline
$K11a364$&$8\,3$&$\{\{11\},\{10,14,16,18,20,22,2,4,6,8,12\}\}$&$\{3,-2,1,-4,5,-6,7,-8,9,-10,11\}$\\ \hline
$K11a365$&$3\,5\,3$&$\{\{11\},\{10,14,16,18,20,22,2,4,8,6,12\}\}$&$\{3,-2,1,-4,7,-6,5,-8,9,-10,11\}$\\ \hline
$K11a366$&$3,3,3++$&$\{\{11\},\{10,14,16,18,20,22,4,2,8,6,12\}\}$&$\{1,-4,3,-2,5,-8,7,-6,11,-10,9\}$\\ \hline
$K11a367$&$11$&$\{\{11\},\{12,14,16,18,20,22,2,4,6,8,10\}\}$&$\{1,-2,3,-4,5,-6,7,-8,9,-10,11\}\}\}$\\ \hline
\end{tabular}

\bigskip

\normalsize

We mentioned that all alternating meander knots are positive, but they represent a proper subset of all positive knots. E.g., the alternating positive knots $K11a43=2\,1,2\,1,2\,1,2$, $K11a123=6^*2.2\,2\,0:2\,0$, $K11a124=6^*2\,1.2:2\,1\,0$, $K11a227=6^*2.2\,1.2:2\,0$, and $K11a329=8^*2:2:.2\,0$ are not meander knots.

The following question can be posed about non-alternating knots with an odd number of crossings: which non-alternating knots have minimal meander diagrams. For $9\le n\le 11$ crossings their list is given in the following table, where every knot is given by its standard symbol or the "Knotscape" ordering number, Conway symbol,  meander Dowker-Thistlethwaite code and short Gauss code (i.e., its meandering part).

\bigskip

\begin{landscape}

\tiny

\noindent \begin{tabular}{|c|c|c|c|}   \hline
$9_{35}$&$3,3,-3$&$\{\{9\},\{8,-12,16,14,18,-4,-2,6,10\}\}$&$\{3,-2,1,-6,5,-4,-9,8,-7\}\}$\\ \hline
$K11n1$&$2\,3,2\,1\,1,-2$&$\{\{11\},\{6,-12,16,22,18,-4,20,2,10,8,14\}\}$&$\{5,-4,3,-8,9,-2,1,-10,-7,6,11\}\}$\\ \hline
$K11n8$&$(2\,2,-2)\,(2\,1,2)$&$\{\{11\},\{6,12,16,22,-18,2,20,4,-8,14,10\}\}$&$\{1,-2,5,-6,11,-10,7,-4,3,8,-9\}\}$\\ \hline
$K11n10$&$(2\,1\,1,-2)\,(3,2)$&$\{\{11\},\{6,-12,16,22,18,-4,20,2,8,10,14\}\}$&$\{1,-8,7,-2,3,-4,-11,10,5,-6,9\}\}$\\ \hline
$K11n11$&$(2\,2,-2)\,(3,2)$&$\{\{11\},\{6,-12,16,22,-18,-4,-20,2,-8,-10,-14\}\}$&$\{1,-8,7,2,-3,4,-11,10,-5,6,9\}\}$\\ \hline
$K11n12$&$2\,3,2\,2,-2$&$\{\{11\},\{6,-12,16,22,-18,-4,-20,2,-10,-8,-14\}\}$&$\{5,-4,3,-8,9,2,-1,10,7,-6,-11\}\}$\\ \hline
$K11n14$&$4\,1,2\,1\,1,-2$&$\{\{11\},\{-6,10,-16,-22,4,-18,-20,-2,-8,-12,-14\}\}$&$\{3,-6,7,-8,9,2,-1,-10,5,-4,11\}\}$\\ \hline
$K11n15$&$4\,1,2\,2,-2$&$\{\{11\},\{-6,10,-16,-22,4,18,20,-2,8,12,14\}\}$&$\{3,6,-7,8,-9,2,-1,10,5,-4,11\}\}$\\ \hline
$K11n17$&$3\,1\,1,2\,1\,1,-2$&$\{\{11\},\{6,-10,-16,22,-4,-18,-20,-2,-8,-14,-12\}\}$&$\{3,-6,9,-8,7,2,-1,-10,5,-4,11\}\}$\\ \hline
$K11n18$&$3\,1\,1,2\,2,-2$&$\{\{11\},\{6,-10,-16,22,-4,18,20,2,-8,14,12\}\}$&$\{3,6,-9,8,-7,2,-1,10,5,-4,11\}\}$\\ \hline
$K11n19$&$5,2\,2,-2$&$\{\{11\},\{-6,10,-16,-22,4,18,20,2,-8,12,14\}\}$&$\{-3,6,-7,8,-9,2,-1,-10,5,-4,11\}\}$\\ \hline
$K11n20$&$3\,2,2\,2,-2$&$\{\{11\},\{6,-10,-16,22,-4,18,20,-2,8,14,12\}\}$&$\{-3,6,-9,8,-7,2,-1,-10,5,-4,11\}\}$\\ \hline
$K11n24$&$(2\,1,2+)\,(3,-2)$&$\{\{11\},\{4,14,18,22,-16,-20,2,-8,-12,6,-10\}\}$&$\{-1,2,-11,8,5,-4,-9,10,3,-6,7\}\}$\\ \hline
$K11n25$&$(3,2+)\,(2\,1,-2)$&$\{\{11\},\{6,-10,18,2,-16,20,22,-4,-8,14,12\}\}$&$\{5,8,-9,2,-3,4,-1,10,7,-6,11\}\}$\\ \hline
$K11n28$&$2\,4,2\,1,-2$&$\{\{11\},\{6,-10,18,2,-16,22,20,-4,-8,14,12\}\}$&$\{5,8,-9,4,-3,2,-1,10,7,-6,11\}\}$\\ \hline
$K11n37$&$(3,2+)\,(3,-2)$&$\{\{11\},\{-6,10,-18,-2,-16,20,22,4,-8,14,12\}\}$&$\{-5,-8,9,2,-3,4,-1,-10,7,-6,-11\}\}$\\ \hline
$K11n38$&$2\,4,3,-2$&$\{\{11\},\{-8,-12,16,20,-2,18,-22,6,4,10,-14\}\}$&$\{-3,2,-1,-8,7,-4,-11,10,5,-6,-9\}\}$\\ \hline
$K11n48$&$3\,1,2\,2,-2\,-1$&$\{\{11\},\{6,-12,16,22,18,-4,-20,-2,8,10,-14\}\}$&$\{-1,8,-7,2,-3,4,11,-10,5,-6,9\}\}$\\ \hline
$K11n49$&$4,2\,2,-2\,-1$&$\{\{11\},\{6,-12,16,22,18,-4,-20,-2,10,8,-14\}\}$&$\{-5,4,-3,-8,9,2,-1,-10,7,-6,-11\}\}$\\ \hline
$K11n50$&$2\,2,2\,1\,1,-2\,-1$&$\{\{11\},\{-4,-12,-22,16,-20,-2,18,6,14,-10,-8\}\}$&$\{3,-4,11,-10,7,-6,5,2,-1,8,-9\}\}$\\ \hline
$K11n51$&$2\,1\,3,2\,1,-2$&$\{\{11\},\{-6,-10,-18,-2,-16,20,22,-4,-8,14,12\}\}$&$\{5,-8,9,2,-3,4,-1,-10,7,-6,11\}\}$\\ \hline
$K11n58$&$3\,1\,2,2\,1,-2$&$\{\{11\},\{-6,-10,18,-2,-16,20,22,4,-8,14,12\}\}$&$\{5,-8,9,2,-3,4,-1,10,7,-6,-11\}\}$\\ \hline
$K11n59$&$3\,1\,1\,1,2\,1,-2$&$\{\{11\},\{6,10,18,2,-16,22,20,4,-8,12,14\}\}$&$\{5,-8,9,-4,1,-2,3,-10,-7,6,11\}\}$\\ \hline
$K11n62$&$4\,2,2\,1,-2$&$\{\{11\},\{-6,-10,18,-2,-16,22,20,4,-8,14,12\}\}$&$\{5,-8,9,4,-3,2,-1,10,7,-6,-11\}\}$\\ \hline
$K11n63$&$4\,1\,1,2\,1,-2$&$\{\{11\},\{6,10,18,2,-16,22,20,4,-8,14,12\}\}$&$\{5,-8,9,-4,3,-2,1,-10,-7,6,11\}\}$\\ \hline
$K11n66$&$(2\,1,-3)\,(2\,1,2)$&$\{\{11\},\{-8,-12,16,18,-22,-4,-20,-2,6,-10,-14\}\}$&$\{-1,2,-3,-6,9,-10,5,-4,11,-8,7\}\}$\\ \hline
$K11n68$&$(3,2\,1)\,(2\,1,-2)$&$\{\{11\},\{-6,-10,-18,-2,-16,22,20,-4,-8,14,12\}\}$&$\{5,-8,9,4,-3,2,-1,-10,7,-6,11\}\}$\\ \hline
$K11n70$&$3\,1\,2,3,-2$&$\{\{11\},\{6,10,-18,2,-16,20,22,-4,-8,14,12\}\}$&$\{-5,8,-9,2,-3,4,-1,-10,7,-6,11\}\}$\\ \hline
$K11n73$&$2\,1,2\,1,-2\,-1,-2$&$\{\{11\},\{6,8,-16,2,20,22,-18,-4,-12,-14,10\}\}$&$\{7,-8,9,4,-5,6,-1,2,-3,-10,11\}\}$\\ \hline
$K11n74$&$2\,1,-2\,-1,2\,1,-2$&$\{\{11\},\{6,8,-16,2,-20,-22,18,-4,12,14,-10\}\}$&$\{-7,8,-9,-4,5,-6,-1,2,-3,-10,11\}\}$\\ \hline
$K11n76$&$3,3,2\,1,-2$&$\{\{11\},\{6,8,-16,2,-20,-22,-18,-4,-12,-14,-10\}\}$&$\{7,-8,9,-4,5,-6,-1,2,-3,-10,11\}\}$\\ \hline
$K11n78$&$3,2\,1,3,-2$&$\{\{11\},\{6,8,16,2,-20,-22,18,4,12,14,-10\}\}$&$\{7,-8,9,4,-5,6,1,-2,3,-10,11\}\}$\\ \hline
$K11n79$&$4\,2,3,-2$&$\{\{11\},\{6,10,-18,2,-16,22,20,-4,-8,14,12\}\}$&$\{-5,8,-9,4,-3,2,-1,-10,7,-6,11\}\}$\\ \hline
$K11n81$&$3,3,3,-2$&$\{\{11\},\{6,8,-16,2,20,22,18,-4,12,14,10\}\}$&$\{7,-8,9,-4,5,-6,1,-2,3,10,-11\}\}$\\ \hline
$K11n82$&$4,2\,2,-3$&$\{\{11\},\{-6,12,16,-22,18,4,-20,-2,8,10,-14\}\}$&$\{-1,-8,7,2,-3,4,-11,10,5,-6,9\}\}$\\ \hline
$K11n83$&$3\,1,2\,2,-3$&$\{\{11\},\{6,-12,-16,22,-18,-4,20,2,-10,-8,14\}\}$&$\{5,-4,3,8,-9,2,-1,10,7,-6,11\}\}$\\ \hline
$K11n84$&$2\,2,2\,2,-3$&$\{\{11\},\{4,12,22,16,-20,2,18,6,14,-10,-8\}\}$&$\{3,-4,11,-10,-7,6,-5,-2,1,-8,9\}\}$\\ \hline
$K11n86$&$-2\,-1\,-1:2:2\,0$&$\{\{11\},\{4,14,18,-22,-16,-20,2,-8,12,6,-10\}\}$&$\{-1,2,-11,-8,5,-4,-9,10,-3,-6,7\}\}$\\ \hline
$K11n87$&$2\,1\,2,2\,1,-3$&$\{\{11\},\{4,14,-18,-22,-16,-20,2,-8,-12,-6,-10\}\}$&$\{-1,2,-11,-8,5,-4,9,-10,3,-6,7\}\}$\\ \hline
$K11n89$&$3\,1,2\,1\,1,-3$&$\{\{11\},\{8,-12,16,18,22,-4,20,2,6,10,14\}\}$&$\{1,-2,3,-6,9,-10,-5,4,11,-8,7\}\}$\\ \hline
$K11n91$&$4,2\,1\,1,-3$&$\{\{11\},\{8,12,16,20,2,18,-22,6,4,10,-14\}\}$&$\{-3,2,-1,-8,7,4,-11,10,-5,6,-9\}\}$\\ \hline
\end{tabular}

\noindent \begin{tabular}{|c|c|c|c|}   \hline
$K11n92$&$-2\,-1\,-1:2:2$&$\{\{11\},\{-8,12,16,18,-22,4,-20,-2,6,10,-14\}\}$&$\{-1,2,-3,-6,9,-10,-5,4,-11,-8,7\}\}$\\ \hline
$K11n96$&$-2\,-1\,-1\,0:2\,0:2\,0$&$\{\{11\},\{-6,-12,16,-22,18,-2,20,4,8,14,-10\}\}$&$\{-1,2,5,-6,11,-10,7,4,-3,8,-9\}\}$\\ \hline
$K11n98$&$(2\,1,-3)\,(3,2)$&$\{\{11\},\{6,10,16,22,4,-18,-20,2,8,-12,-14\}\}$&$\{3,6,-7,8,-9,-2,1,-10,5,-4,11\}\}$\\ \hline
$K11n99$&$(3,-2\,-1)\,(3,2)$&$\{\{11\},\{-8,-12,16,-20,-2,-18,-22,6,4,-10,-14\}\}$&$\{3,-2,1,8,-7,4,-11,10,-5,6,-9\}\}$\\ \hline
$K11n101$&$2\,3,2\,1,-3$&$\{\{11\},\{6,10,16,22,4,-18,-20,2,8,-14,-12\}\}$&$\{3,6,-9,8,-7,-2,1,-10,5,-4,11\}\}$\\ \hline
$K11n103$&$2\,1\,1,2\,1\,1,-3$&$\{\{11\},\{4,14,-18,22,16,20,2,8,12,-6,10\}\}$&$\{1,-2,11,-8,5,-4,-9,10,3,-6,7\}\}$\\ \hline
$K11n105$&$2\,1\,1\,1,2\,1,-3$&$\{\{11\},\{4,14,18,22,16,-20,2,8,12,6,-10\}\}$&$\{1,-2,11,-8,-5,4,9,-10,3,-6,7\}\}$\\ \hline
$K11n108$&$-3\,0:2\,1:2$&$\{\{11\},\{-6,12,16,-22,18,2,20,4,8,14,10\}\}$&$\{1,-2,5,-6,-11,10,7,-4,3,-8,9\}\}$\\ \hline
$K11n109$&$-2\,-1\,0:-3\,0:-2\,0$&$\{\{11\},\{-6,12,16,-22,18,4,20,2,8,10,14\}\}$&$\{1,8,-7,-2,3,-4,11,-10,5,-6,9\}\}$\\ \hline
$K11n110$&$-2\,-1\,-1\,0:2:2\,0$&$\{\{11\},\{6,-12,16,22,18,-2,-20,4,8,-14,-10\}\}$&$\{-1,2,5,-6,-11,10,7,-4,3,8,-9\}\}$\\ \hline
$K11n111$&$-3\,-1\,0:2\,0:2\,0$&$\{\{11\},\{-6,12,16,-22,-18,4,-20,2,-8,-10,-14\}\}$&$\{1,8,-7,2,-3,4,11,-10,-5,6,9\}\}$\\ \hline
$K11n113$&$-2\,-2:2:2$&$\{\{11\},\{6,-12,-16,22,-18,-4,-20,-2,-10,-8,-14\}\}$&$\{5,-4,3,-8,9,2,-1,-10,7,-6,11\}\}$\\ \hline
$K11n114$&$-4:2:2\,0$&$\{\{11\},\{-8,12,16,-20,-2,18,-22,6,4,10,-14\}\}$&$\{-3,2,-1,-8,7,4,-11,10,5,-6,9\}\}$\\ \hline
$K11n116$&$-3\,0:2\,1:-2\,0$&$\{\{11\},\{6,-12,-16,22,18,-4,-20,-2,10,8,-14\}\}$&$\{-5,4,-3,-8,9,2,-1,-10,7,-6,11\}\}$\\ \hline
$K11n117$&$-4\,0:2\,0:2\,0$&$\{\{11\},\{6,-12,-16,22,18,-4,20,-2,10,8,14\}\}$&$\{-5,4,-3,8,-9,2,-1,-10,7,-6,11\}\}$\\ \hline
$K11n121$&$-2\,-2:-2\,0:-2\,0$&$\{\{11\},\{6,12,16,22,18,4,-20,2,8,10,-14\}\}$&$\{1,-8,7,-2,3,-4,11,-10,-5,6,9\}\}$\\ \hline
$K11n122$&$3\,2,3,-3$&$\{\{11\},\{8,14,-18,-16,22,20,4,2,-6,10,12\}\}$&$\{5,-4,1,-2,3,-8,7,-6,-11,10,-9\}\}$\\ \hline
$K11n131$&$-3\,0:2\,1\,0:2$&$\{\{11\},\{6,12,16,22,-18,4,20,2,-8,-10,14\}\}$&$\{1,-8,7,2,-3,4,11,-10,5,-6,9\}\}$\\ \hline
$K11n134$&$2.-2\,-1.2.2$&$\{\{11\},\{6,12,16,22,-18,4,20,2,-10,-8,14\}\}$&$\{-5,4,-3,-8,9,-2,1,-10,7,-6,11\}\}$\\ \hline
$K11n137$&$3\,1\,1,2\,1,-3$&$\{\{11\},\{-6,10,18,-2,16,22,20,4,8,12,14\}\}$&$\{5,8,-9,-4,1,-2,3,-10,7,-6,11\}\}$\\ \hline
$K11n138$&$3\,1\,1,3,-2\,-1$&$\{\{11\},\{8,-14,-18,-16,20,22,-4,-2,-6,12,10\}\}$&$\{-5,2,-3,4,-1,-8,7,-6,11,-10,9\}\}$\\ \hline
$K11n139$&$5,3,-3$&$\{\{11\},\{8,-14,18,16,22,20,-4,-2,6,12,10\}\}$&$\{5,-4,3,-2,1,-8,7,-6,-11,10,-9\}\}$\\ \hline
$K11n140$&$4\,1,2\,1,-3$&$\{\{11\},\{-8,12,16,20,-2,18,22,6,4,10,14\}\}$&$\{-3,2,-1,8,-7,4,-11,10,5,-6,-9\}\}$\\ \hline
$K11n141$&$4\,1,3,-2\,-1$&$\{\{11\},\{8,-14,-18,-16,22,20,-4,-2,-6,12,10\}\}$&$\{-5,4,-3,2,-1,-8,7,-6,11,-10,9\}\}$\\ \hline
$K11n143$&$2.-2\,-1.-2\,0.2\,0$&$\{\{11\},\{-6,-12,16,-22,18,-2,-20,4,8,-14,-10\}\}$&$\{-1,2,5,-6,11,-10,7,-4,3,8,-9\}\}$\\ \hline
$K11n145$&$2.-3.-2\,0.2\,0$&$\{\{11\},\{6,12,16,22,-18,4,-20,2,-8,-10,-14\}\}$&$\{-1,6,-5,2,-11,10,3,-4,7,-8,9\}\}$\\ \hline
$K11n164$&$2.2.-2.2\,0.2\,0$&$\{\{11\},\{-6,-12,16,-22,-18,-2,-20,4,-8,-14,-10\}\}$&$\{-1,2,5,-6,11,-10,7,-4,3,-8,9\}\}$\\ \hline
$K11n170$&$2\,0.3.-2.2$&$\{\{11\},\{6,12,16,22,18,4,-20,2,10,8,-14\}\}$&$\{1,-6,5,-4,-9,10,3,-2,11,-8,7\}$\\ \hline
\end{tabular}

\end{landscape}

\normalsize

\begin{figure}[th]
\centerline{\psfig{file=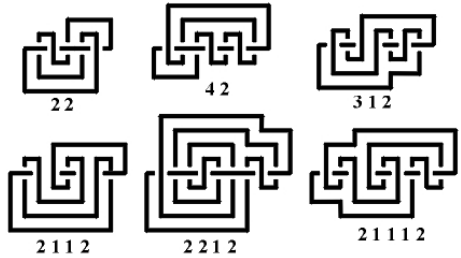,width=2.4in}} \vspace*{8pt}
\caption{Non-minimal meander diagrams of knots $4_1$, $6_1$, $6_2$, $6_3$, $7_6$, and $7_7$.
\label{f8}}
\end{figure}

The next question is: which knots can be represented by non-minimal meander diagrams. E.g., figure-eight knot $4_1=2\,2$ can be represented by non-minimal meander diagram given by Gauss code $\{-1,2,-3,-4,5,3,-2,1,4,-5\}$ with $n=5$ crossings (Fig. \ref{f8}). This means that for every knot which is not meander knot (i.e., knot with a minimal meander diagram), and which can be represented by some non-minimal meander diagram we can define its {\it meander number}: the minimum number of crossings of its meander diagrams, where the minimum is taken over all its meander diagrams. However, we first need to check which knots have meander diagrams. Except alternating meander knots, we are sure that non-alternating knots with the same shadows will have the meander diagram, with some crossings changed from overcrossings to undercrossings or {\it vice versa}. The next step is to make all possible crossing changes in alternating minimal meander diagrams, i.e., in Gauss codes of alternating meander knots an see which knots will be obtained.

\begin{theorem}
Every knot with at most $n=9$ crossings has a meander diagram.
\end{theorem}

\noindent {\bf Question:} Is it true that every (alternating) knot has a meander diagram?

\bigskip

For $n=10$, we have not succeeded to find meander diagrams only for three knots: $10_{96}=.2.2\,1.2$, $10_{99}=.2.2.2\,0.2\,0$, and $10_{123}=10^*$.

In the following table are given non-alternating meander diagrams of all knots with $n\le 9$ crossings which are not alternating meander knots. As a representative of meander diagrams, for every knot is taken a non-alternating meander diagram with the minimal number of crossings. Every knot is given by its standard symbol, Conway symbol, meander Dowker-Thisltlethwaite code, short Gauss code, and meander number. For non-alternating meander diagrams, short Gausse code is the meander part of the Gauss code, from which the complete Gauss code can be reconstructed: the first part of the complete code will be the sequence of the numbers from $1$ to $n$, where each of them has the opposite sign then the same number from the given meander sequence. E.g., short Gauss code of the figure-eight knot is $\{3,-2,1,4,-5\}$, and its complete Gauss code will be $\{-1,2,-3,-4,5,3,-2,1,4,-5\}$.

\begin{landscape}

\tiny

\noindent \begin{tabular}{|c|c|c|c|c|}   \hline
$4_1$&$2\,2$&$\{\{5\},\{4,-8,10,-2,6\}\}$&$\{3,-2,1,4,-5\}$&5\\ \hline
$6_1$&$4\,2$&$\{\{7\},\{-4,10,-14,-12,2,-8,-6\}\}$&$\{-1,2,7,-6,5,-4,3\}$&7\\ \hline
$6_2$&$3\,1\,2$&$\{\{7\},\{6,-10,-12,14,-2,-4,8\}\}$&$\{-3,2,-1,-4,5,-6,7\}$&7\\ \hline
$6_3$&$2\,1\,1\,2$&$\{\{7\},\{4,10,14,-12,2,-6,-8\}\}$&$\{-3,4,-5,-2,1,-6,7\}$&7\\ \hline
$7_6$&$2\,2\,1\,2$&$\{\{9\},\{6,-12,14,-16,-18,-4,-2,10,-8\}\}$&$\{-1,6,3,-4,5,2,9,-8,7\}$&9\\ \hline
$7_7$&$2\,1\,1\,1\,2$&$\{\{9\},\{-8,-12,14,16,-18,-4,-2,6,-10\}\}$&$\{-1,2,-3,-6,5,-4,9,-8,7\}$&9\\ \hline
$8_1$&$6\,2$&$\{\{9\},\{-4,12,-18,-16,-14,2,-10,-8,-6\}\}$&$\{-1,2,9,-8,7,-6,5,-4,3\}$&9\\ \hline
$8_2$&$5\,1\,2$&$\{\{9\},\{8,-12,-14,-16,18,-2,-4,-6,10\}\}$&$\{-3,2,-1,-4,5,-6,7,-8,9\}$&9\\ \hline
$8_3$&$4\,4$&$\{\{9\},\{-6,12,14,-18,-16,4,2,-10,-8\}\}$&$\{-1,4,-3,2,9,-8,7,-6,5\}$&9\\ \hline
$8_4$&$4\,1\,3$&$\{\{9\},\{-6,12,14,-16,-18,4,2,-10,-8\}\}$&$\{-1,-6,3,-4,5,-2,-9,8,-7\}$&9\\ \hline
$8_5$&$3,3,2$&$\{\{9\},\{4,-12,16,18,14,-2,8,10,6\}\}$&$\{5,-6,7,-2,3,-4,1,8,-9\}$&9\\ \hline
$8_6$&$3\,3\,2$&$\{\{9\},\{6,-12,14,16,18,-4,-2,10,8\}\}$&$\{1,-6,3,-4,5,-2,-9,8,-7\}$&9\\ \hline
$8_7$&$4\,1\,1\,2$&$\{\{9\},\{4,12,18,-14,-16,2,-6,-8,-10\}\}$&$\{3,-4,5,-6,7,2,-1,8,-9\}$&9\\ \hline
$8_8$&$2\,3\,1\,2$&$\{\{9\},\{4,10,18,-14,2,-16,-6,-12,-8\}\}$&$\{7,-6,-1,2,5,-8,9,4,-3\}$&9\\ \hline
$8_9$&$3\,1\,1\,3$&$\{\{9\},\{6,-12,-14,18,16,-2,-4,8,10\}\}$&$\{3,-4,5,-2,1,6,-7,8,-9\}$&9\\ \hline
$8_{10}$&$2\,1,3,2$&$\{\{9\},\{4,12,-16,-18,14,2,8,10,-6\}\}$&$\{5,-6,7,2,-3,4,1,-8,9\}$&9\\ \hline
$8_{11}$&$3\,2\,1\,2$&$\{\{9\},\{-8,12,14,16,-18,4,2,6,-10\}\}$&$\{1,-2,3,6,-5,4,9,-8,7\}$&9\\ \hline
$8_{12}$&$2\,2\,2\,2$&$\{\{11\},\{-4,12,-22,-18,-16,2,20,6,-8,14,-10\}\}$&$\{9,-8,-1,-4,3,-2,7,10,-11,6,-5\}$&11\\ \hline
$8_{13}$&$3\,1\,1\,1\,2$&$\{\{9\},\{6,12,14,-16,-18,4,2,10,-8\}\}$&$\{1,-6,-3,4,-5,-2,9,-8,7\}$&9\\ \hline
$8_{14}$&$2\,2\,1\,1\,2$&$\{\{9\},\{4,12,18,14,-16,2,6,-10,-8\}\}$&$\{-5,4,-3,-6,7,-2,1,-8,9\}$&9\\ \hline
$8_{15}$&$2\,1,2\,1,2$&$\{\{11\},\{4,-14,18,22,16,20,2,8,12,6,10\}\}$&$\{7,-4,3,-8,9,-2,5,-6,1,10,11\}$&11\\ \hline
$8_{16}$&$.2.2\,0$&$\{\{11\},\{-6,10,-18,-2,16,-20,-22,4,8,14,-12\}\}$&$\{-1,8,5,-6,7,2,-11,10,3,-4,9\}$&11\\ \hline
$8_{17}$&$.2.2$&$\{\{11\},\{4,14,-18,-22,-16,20,2,-8,12,-6,10\}\}$&$\{7,-4,3,8,-9,2,-5,6,1,-10,11\}$&11\\ \hline
$8_{18}$&$8^*$&$\{\{11\},\{-6,12,-16,-22,18,2,-20,-4,8,-14,10\}\}$&$\{7,-6,-1,2,-5,8,-9,-4,3,-10,11\}$&11\\ \hline
$8_{19}$&$3,3,-2$&$\{\{9\},\{4,-12,-16,-18,-14,-2,-8,-10,-6\}\}$&$\{5,-6,7,-2,3,-4,-1,-8,9\}$&9\\ \hline
$8_{20}$&$3,2\,1,-2$&$\{\{9\},\{4,-12,16,18,-14,-2,-8,-10,6\}\}$&$\{5,-6,7,2,-3,4,-1,-8,9\}$&9\\ \hline
$8_{21}$&$2\,1,2\,1,-2$&$\{\{9\},\{4,12,-16,-18,-14,2,-8,-10,-6\}\}$&$\{5,-6,7,-2,3,-4,-1,8,-9\}$&9\\ \hline
$9_8$&$2\,4\,1\,2$&$\{\{11\},\{-4,-12,-22,-18,-16,-2,20,6,-8,14,10\}\}$&$\{9,-8,-1,-4,3,-2,-7,-10,11,6,-5\}$&11\\ \hline
$9_{11}$&$4\,1\,2\,2$&$\{\{11\},\{-4,12,-22,16,18,2,20,6,8,14,-10\}\}$&$\{7,-6,1,-2,-5,-8,9,-10,11,4,-3\}$&11\\ \hline
$9_{12}$&$4\,2\,1\,2$&$\{\{11\},\{4,12,22,18,16,2,-20,6,8,-14,-10\}\}$&$\{9,-8,1,-4,3,-2,-7,-10,11,6,-5\}$&11\\ \hline
$9_{14}$&$4\,1\,1\,1\,2$&$\{\{11\},\{-4,12,-22,18,16,2,-20,6,8,-14,-10\}\}$&$\{-9,8,1,-4,3,-2,-7,-10,11,6,-5\}$&11\\ \hline
$9_{15}$&$2\,3\,2\,2$&$\{\{11\},\{-4,-12,-22,-18,-16,-2,20,6,-8,14,-10\}\}$&$\{9,-8,-1,-4,3,-2,7,-10,11,6,-5\}$&11\\ \hline
$9_{17}$&$2\,1\,3\,1\,2$&$\{\{11\},\{-4,12,-22,16,-20,2,-18,6,-14,-8,-10\}\}$&$\{9,-8,-1,2,5,-6,7,10,-11,-4,3\}$&11\\ \hline
$9_{19}$&$2\,3\,1\,1\,2$&$\{\{11\},\{4,-12,22,-18,-16,-2,20,6,-8,14,10\}\}$&$\{9,-8,1,4,-3,2,7,10,-11,-6,5\}$&11\\ \hline
$9_{20}$&$3\,1\,2\,1\,2$&$\{\{11\},\{4,12,22,16,18,2,-20,6,8,-14,-10\}\}$&$\{-7,6,1,-2,-5,-8,9,-10,11,-4,3\}$&11\\ \hline
$9_{21}$&$3\,1\,1\,2\,2$&$\{\{11\},\{-4,12,-22,18,16,2,20,6,8,14,-10\}\}$&$\{-9,8,1,-4,3,-2,-7,-10,11,-6,5\}$&11\\ \hline
$9_{22}$&$2\,1\,1,3,2$&$\{\{11\},\{4,-14,20,22,16,-18,-2,8,-12,-10,6\}\}$&$\{-7,6,-5,-8,9,-2,3,-4,1,10,-11\}$&11\\ \hline
$9_{24}$&$2\,1,3,2+$&$\{\{11\},\{-4,14,-22,-18,-20,16,2,10,12,-8,-6\}\}$&$\{-1,-10,9,-8,5,-6,7,2,-3,4,-11\}$&11\\ \hline
$9_{25}$&$2\,2,2\,1,2$&$\{\{11\},\{-6,12,16,-22,-18,4,-20,-2,-10,-8,-14\}\}$&$\{-1,6,-5,-2,11,-10,3,-4,9,-8,7\}$&11\\ \hline
$9_{26}$&$3\,1\,1\,1\,1\,2$&$\{\{11\},\{-4,12,-22,16,18,2,-20,6,8,-14,-10\}\}$&$\{-7,6,1,-2,-5,-8,9,-10,11,4,-3\}$&11\\ \hline
$9_{27}$&$2\,1\,2\,1\,1\,2$&$\{\{11\},\{-4,12,-22,16,-20,2,18,6,14,-8,-10\}\}$&$\{9,-8,1,-2,-5,6,-7,-10,11,4,-3\}$&11\\ \hline
$9_{28}$&$2\,1,2\,1,2+$&$\{\{11\},\{4,14,22,-18,-20,-16,2,-10,-12,-8,6\}\}$&$\{-1,10,-9,8,5,-6,7,-2,3,-4,-11\}$&11\\ \hline
$9_{29}$&$.2.2\,0.2$&$\{\{13\},\{-4,16,20,-22,-26,-24,18,2,12,14,-8,-6,-10\}\}$&$\{-1,2,-11,12,-13,-6,5,-4,9,-8,7,10,3\}$&13\\ \hline
$9_{30}$&$2\,1\,1,2\,1,2$&$\{\{11\},\{4,14,-20,-22,16,-18,2,8,-12,-10,-6\}\}$&$\{7,-6,5,8,-9,-2,3,-4,-1,10,-11\}$&11\\ \hline
\end{tabular}

\noindent \begin{tabular}{|c|c|c|c|c|}   \hline
$9_{31}$&$2\,1\,1\,1\,1\,1\,2$&$\{\{11\},\{4,12,22,16,-20,2,18,6,14,-8,-10\}\}$&$\{9,-8,1,-2,-5,6,-7,-10,11,-4,3\}$&11\\ \hline
$9_{32}$&$.2\,1.2\,0$&$\{\{11\},\{-6,12,-16,-22,-18,2,20,-4,-8,14,-10\}\}$&$\{-7,6,-1,2,5,-8,9,-4,3,-10,11\}$&11\\ \hline
$9_{33}$&$.2\,1.2$&$\{\{11\},\{6,-12,16,22,-18,-4,20,2,-8,-10,14\}\}$&$\{1,-8,7,2,-3,4,-11,10,5,-6,9\}$&11\\ \hline
$9_{34}$&$8^*2\,0$&$\{\{13\},\{-6,-12,18,-26,22,-4,20,-24,2,-10,14,8,-16\}\}$&$\{-1,-8,7,2,13,-12,3,-4,-11,10,-5,6,9\}$&13\\ \hline
$9_{36}$&$2\,2,3,2$&$\{\{11\},\{4,14,20,22,-16,18,2,-8,12,10,6\}\}$&$\{7,-6,5,8,-9,-2,3,-4,1,-10,11\}$&11\\ \hline
$9_{37}$&$2\,1,2\,1,3$&$\{\{11\},\{-4,12,-22,16,-20,2,-18,6,-14,-10,-8\}\}$&$\{9,-8,-1,2,7,-6,5,10,-11,-4,3\}$&11\\ \hline
$9_{39}$&$2:2:2\,0$&$\{\{13\},\{-8,-18,-14,20,-26,-22,-6,24,-4,-2,-12,-10,16\}\}$&$\{-1,-8,7,-4,3,-2,13,-12,-5,6,11,-10,9\}$&13\\ \hline
$9_{40}$&$9^*$&$\{\{13\},\{-6,-14,18,-26,-20,24,-4,-22,2,-8,12,-16,10\}\}$&$\{-1,-10,9,-2,3,8,-7,4,13,-12,5,-6,-11\}$&13\\ \hline
$9_{41}$&$2\,0:2\,0:2\,0$&$\{\{11\},\{-6,12,16,-22,18,4,-20,2,10,8,-14\}\}$&$\{1,-6,5,-2,-11,10,-3,4,9,-8,7\}$&11\\ \hline
$9_{42}$&$2\,2,3,-2$&$\{\{13\},\{-6,-10,20,-2,18,26,-24,-22,4,8,-16,-14,12\}\}$&$\{-1,6,-5,2,13,-12,11,-10,-9,8,3,-4,7\}$&13\\ \hline
$9_{43}$&$2\,1\,1,3,-2$&$\{\{11\},\{4,-14,-20,-22,-16,18,-2,-8,12,10,-6\}\}$&$\{-7,6,-5,-8,9,-2,3,-4,-1,-10,11\}$&11\\ \hline
$9_{45}$&$2\,1\,1,2\,1,-2$&$\{\{13\},\{6,10,20,2,-18,-26,24,22,4,-8,16,14,-12\}\}$&$\{1,6,-5,-2,13,-12,11,-10,-9,8,3,-4,7\}$&13\\ \hline
$9_{46}$&$3,3,-3$&$\{\{11\},\{-4,12,-22,-16,-20,2,18,-6,14,-10,-8\}\}$&$\{-9,8,-1,2,7,-6,5,-10,11,-4,3\}$&11\\ \hline
$9_{47}$&$8^*-2\,0$&$\{\{11\},\{-6,12,16,-22,18,2,-20,4,8,-14,-10\}\}$&$\{-7,6,1,-2,-5,-8,9,4,-3,-10,11\}$&11\\ \hline
$9_{48}$&$2\,1,2\,1,-3$&$\{\{13\},\{-6,10,20,-2,18,-26,24,22,4,8,16,14,-12\}\}$&$\{1,-6,5,-2,13,-12,11,-10,-9,8,-3,4,7\}$&13\\ \hline
$9_{49}$&$-2\,0:-2\,0:-2\,0$&$\{\{11\},\{4,-14,-18,22,-16,-20,-2,-8,-12,-6,-10\}\}$&$\{7,-4,3,-8,9,2,5,-6,-1,-10,11\}$&11\\ \hline
\end{tabular}

\end{landscape}

\normalsize

\section{2-component meander links}

Open meanders with an even number of crossings offer another interesting possibility: to join in pairs loose ends of the meander axis, and loose ends of the meander curve. As a result, we obtain a shadow of a 2-component link with one component in the form of a circle, and the other component meandering around it. Natural question is: which alternating links can be obtained from these shadows, and, in general: which 2-component links have meander diagrams. It is clear that every component has no self-intersections, so the set of 2-component meander link coincides with the set of alternating 2-component links with components without self-intersections, and all their minimal diagrams will preserve this property. Among meander links we can also obtain non-prime meander links, so our discussion we will restrict only to prime meander link diagrams. Every minimal diagram of a 2-component meander link is a meander diagram.

\begin{figure}[th]
\centerline{\psfig{file=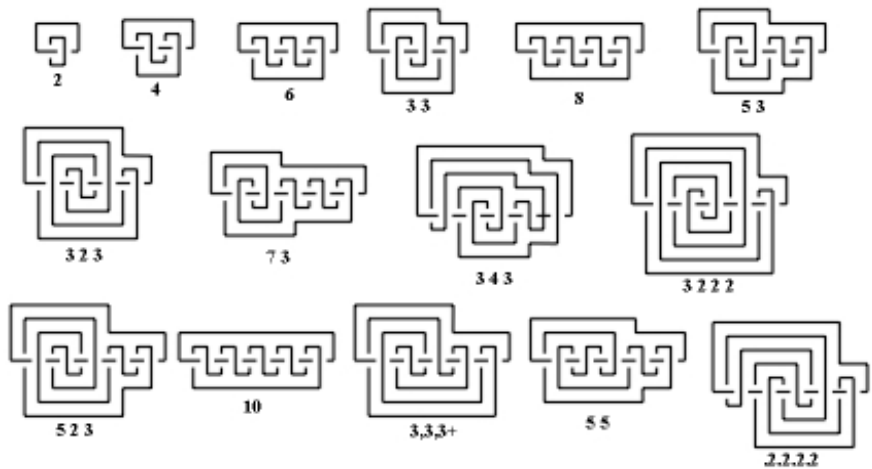,width=4.6in}} \vspace*{8pt}
\caption{Meander links up to $n=10$ crossings.
\label{links01}}
\end{figure}

In the following table is given the number of open meanders $OM$ with $n$ crossings for even $n$, and the number $AML$ of their corresponding alternating prime meander links.

\small

\bigskip

\begin{tabular}{|c|c|c|}   \hline
$n$ & $OM$ & $AML$  \\  \hline
2 & 1 &  1  \\  \hline
4 & 3 &  1  \\  \hline
6 & 14 &   2 \\  \hline
8 & 81 &   3  \\  \hline
10 & 538 &   8  \\  \hline
12 & 3926 &  17  \\  \hline
14& 30694 &  56   \\  \hline
16 & 252939 & 202  \\  \hline
\end{tabular}

\bigskip

\normalsize

In the following table is given the complete list of alternating meander links up to $n=12$ crossings. Meander links are given by their symbols (up to $n=10$ crossings), Dowker-Thisltlethwaite meander codes, and short Gauss codes, where the first sequence $(-1)^ii$ ($1\le i\le n$) corresponding to the axis component is omitted, so every meander link is given by its corresponding meander permutation corresponding to the other component.

\bigskip

\tiny

\noindent \begin{tabular}{|c|c|c|c|}   \hline
$2_1^2$&2&$\{\{1,1\},\{4,2\}\}$&$\{1,-2\}$\\ \hline
$4_1^2$&$4$&$\{\{2,2\},\{6,8,2,4\}\}$&$\{1,-2,3,-4\}$\\ \hline
$6_1^2$&$6$&$\{\{3,3\},\{8,10,12,6,2,4\}\}$&$\{1,-2,3,-4,5,-6\}$\\ \hline
$6_2^2$&$3\,3$&$\{\{3,3\},\{8,10,12,2,6,4\}\}$&$\{1,-4,3,-2,5,-6\}$\\ \hline
$8_1^2$&$8$&$\{\{4,4\},\{10,12,14,16,8,2,4,6\}\}$&$\{1,-2,3,-4,5,-6,7,-8\}$\\ \hline
$8_2^2$&$5\,3$&$\{\{4,4\},\{10,12,14,16,2,8,4,6\}\}$&$\{1,-4,3,-2,5,-6,7,-8\}$\\ \hline
$8_4^2$&$3\,2\,3$&$\{\{4,4\},\{10,12,14,16,2,8,6,4\}\}$&$\{1,-6,3,-4,5,-2,7,-8\}$\\ \hline
$L10a114$&$7\,3$&$\{\{5,5\},\{12,14,16,18,20,2,10,4,6,8\}\}$&$\{1,-4,3,-2,5,-6,7,-8,9,-10\}$\\ \hline
$L10a115$&$3\,4\,3$&$\{\{5,5\},\{12,14,16,18,20,2,10,4,8,6\}\}$&$\{-2,9,-8,7,-4,5,-6,3,-10,1\}$\\ \hline
$L10a116$&$3\,2\,2\,3$&$\{\{5,5\},\{12,14,16,18,20,2,10,8,6,4\}\}$&$\{1,-8,3,-6,5,-4,7,-2,9,-10\}$\\ \hline
$L10a117$&$5\,2\,3$&$\{\{5,5\},\{12,14,16,18,20,4,2,10,6,8\}\}$&$\{1,-6,3,-4,5,-2,7,-8,9,-10\}$\\ \hline
$L10a118$&$10$&$\{\{5,5\},\{12,14,16,18,20,10,2,4,6,8\}\}$&$\{1,-2,3,-4,5,-6,7,-8,9,-10\}$\\ \hline
$L10a119$&$3,3,3+$&$\{\{5,5\},\{12,14,16,20,18,2,10,6,4,8\}\}$&$\{1,-6,7,-8,3,-4,5,-2,9,-10\}$\\ \hline
$L10a120$&$5\,5$&$\{\{5,5\},\{12,14,16,20,18,6,2,4,10,8\}\}$&$\{1,-6,5,-4,3,-2,7,-8,9,-10\}$\\ \hline
$L10a121$&$.2.2.2.2$&$\{\{5,5\},\{12,16,20,14,18,2,6,10,4,8\}\}$&$\{-2,7,-6,3,-10,9,-4,5,-8,1\}$\\ \hline
&$12$&$\{\{6,6\},\{14,16,18,20,22,24,12,2,4,6,8,10\}\}$&$\{1,-2,3,-4,5,-6,7,-8,9,-10,11,-12\}$\\ \hline
&$3\,2\,2\,2\,3$&$\{\{6,6\},\{14,16,18,20,22,24,2,12,10,8,6,4\}\}$&$\{1,-10,3,-8,5,-6,7,-4,9,-2,11,-12\}$\\ \hline
&$3\,2,3,3+$&$\{\{6,6\},\{14,16,18,20,24,22,2,12,8,6,4,10\}\}$&$\{1,-2,3,-10,11,-12,5,-8,7,-6,9,-4\}$\\ \hline
&$3,3,3+++$&$\{\{6,6\},\{14,16,18,20,22,24,2,12,6,4,10,8\}\}$&$\{-2,11,-8,9,-10,7,-4,5,-6,3,-12,1\}$\\ \hline
&$3,3,3,3$&$\{\{6,6\},\{14,16,24,20,22,18,6,2,10,12,8,4\}\}$&$\{1,-8,9,-10,5,-6,7,-2,3,-4,11,-12\}$\\ \hline
&$3\,4\,2\,3$&$\{\{6,6\},\{14,16,18,20,22,24,2,12,4,10,8,6\}\}$&$\{-2,11,-10,9,-4,7,-6,5,-8,3,-12,1\}$\\ \hline
&$3\,6\,3$&$\{\{6,6\},\{14,16,18,20,22,24,2,12,4,6,10,8\}\}$&$\{-2,11,-10,9,-8,7,-4,5,-6,3,-12,1\}$\\ \hline
&$5\,2\,2\,3$&$\{\{6,6\},\{14,16,18,20,22,24,2,12,10,6,8,4\}\}$&$\{1,-8,3,-6,5,-4,7,-2,9,-10,11,-12\}$\\ \hline
&$5\,2\,5$&$\{\{6,6\},\{14,16,18,20,22,24,6,2,4,12,8,10\}\}$&$\{1,-8,3,-4,5,-6,7,-2,9,-10,11,-12\}$\\ \hline
&$5,3,3+$&$\{\{6,6\},\{14,16,18,22,24,20,6,2,4,10,12,8\}\}$&$\{1,-6,7,-8,3,-4,5,-2,9,-10,11,-12\}$\\ \hline
&$5\,4\,3$&$\{\{6,6\},\{14,16,18,20,24,22,8,2,6,4,12,10\}\}$&$\{1,-8,7,-6,3,-4,5,-2,9,-10,11,-12\}$\\ \hline
&$6^*2.2.3.2\,0.2\,0$&$\{\{6,6\},\{14,16,20,24,18,22,6,2,10,4,8,12\}\}$&$\{-2,7,-6,3,-12,11,-4,5,-10,9,-8,1\}$\\ \hline
&$6^*3\,1.2:2\,0.2\,0$&$\{\{6,6\},\{14,16,20,24,18,22,2,12,6,10,4,8\}\}$&$\{-2,9,-6,7,-8,3,-12,11,-4,5,-10,1\}$\\ \hline
&$6^*4.2:2\,0.2\,0$&$\{\{6,6\},\{14,16,20,24,18,22,10,2,4,8,12,6\}\}$&$\{-2,7,-6,3,-12,11,-10,9,-4,5,-8,1\}$\\ \hline
&$7\,2\,3$&$\{\{6,6\},\{14,16,18,20,22,24,4,2,12,6,8,10\}\}$&$\{1,-6,3,-4,5,-2,7,-8,9,-10,11,-12\}$\\ \hline
&$7\,5$&$\{\{6,6\},\{14,16,18,20,24,22,8,2,4,6,12,10\}\}$&$\{1,-6,5,-4,3,-2,7,-8,9,-10,11,-12\}$\\ \hline
&$9\,3$&$\{\{6,6\},\{14,16,18,20,22,24,2,12,4,6,8,10\}\}$&$\{1,-4,3,-2,5,-6,7,-8,9,-10,11,-12\}$\\ \hline
\end{tabular}

\normalsize

\bigskip

Same as for knots, we can pose for 2-component links the natural question: which 2-component links have meander diagrams. From the definition of meander links it is clear that this will be links with both components which are unknots. Moreover, every component will have no self-crossings, i.e., it will be a shadow of a circle. In the case of alternating minimal meander diagrams, all such diagrams of a 2-component link will have this property. However, in the case of non-minimal meander diagrams, some links with an odd number of crossings can be represented by meander diagrams. Moreover, their minimal diagrams have components with self-intersections, but in their non-minimal meander diagrams none of components has self-intersections. E.g., Whitehead link $5_1^2=2\,1\,2$ can be represented by non-minimal meander diagram with the Gauss code
$\{\{1,-2,3,-4,-5,6,-7,8,9,10\},\{-1,-6,7,-8,-3,4,5,2,-9,-10\}$ with $n=10$ crossings. In the same way, by crossing changes, we obtained non-minimal meander diagrams of 2-component links with an odd number of crossings: $7_1^2$, $7_2^2$, $7_3^2$, $7_6^2$, $9_1^2$, $9_2^2$, $9_3^2$, $9_4^2$, $9_5^2$, $9_6^2$, $9_7^2$, $9_8^2$, $9_9^2$, $9_{11}^2$, $9_{12}^2$, $9_{24}^2$, $9_{34}^2$, $9_{35}^2$, $9_{41}^2$, $9_{42}^2$, $9_{53}^2$, $9_{54}^2$, $9_{57}^2$, $9_{58}^2$, and $9_{61}^2$.

\section{Multi-component meander links}

\begin{figure}[th]
\centerline{\psfig{file=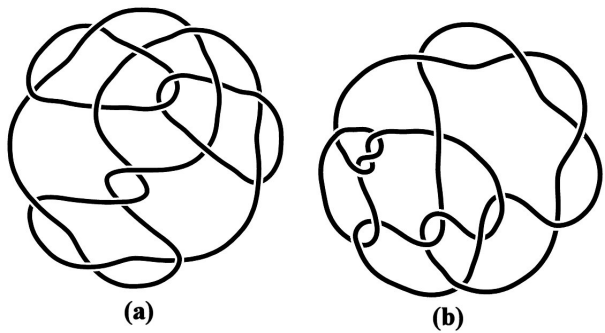,width=3.2in}} \vspace*{8pt}
\caption{A 3-component and 4-component meander link with $n=16$ crossings.
\label{multi}}
\end{figure}

In this section we will consider meander links with at least 3 components. An alternating link will be called meander link if it has an alternating minimal diagram with the property that every pair of its components represents a 2-component meander link (Fig. \ref{multi}). According to this definition every component has no self intersection.  Moreover, in the list of meander links will de not included links with some disjoint components or links where some components make 2-component split link. E.g., in this list will be not included link $8_1^4=2,2,2,2$ since it contains pairs of disjoint components, or Borromean rings $6_2^3=6^*$ since they contain pairs of components making a split link. Every multi-component meander link has an even number of crossings.

For $n\le 16$ crossings we have the following numbers of the $c$-component prime alternating meander links:

\bigskip

\begin{tabular}{|c|c|c|}   \hline
$n$ & $c=3$ & $c=4$  \\  \hline
6 & 1 &    \\  \hline
8 & 2 &    \\  \hline
10 & 12 &   \\  \hline
12 & 59 &  4   \\  \hline
14 & 383 &  32   \\  \hline
16 & 3200 &  388  \\  \hline
\end{tabular}

\bigskip

Multi-component meander links are selected from the unpublished lists of links derived by Morwen Thistlethwaite, using the criteria that every component of a multi-component meander link has no self-intersection and that every 2-component sublink of a meander multi-component link is a 2-component meander link.

In the following tables are given all meander links with at most $n=12$ crossings. Every link is given by its standard symbol from the Rolfsen's book \cite{rol} or from Thistlethwaite's list, by its Conway symbol, meander Dowker-Thistethwaite code and Gauss code.

\bigskip

\tiny

\noindent \begin{tabular}{|c|c|c|c|}   \hline
$6_1^3$&$2,2,2$&$\{\{2,2,2\},\{8,10,2,12,6,4\}\}$&$\{\{1,-3,2,-6\},\{3,-5,4,-1\},\{5,-2,6,-4\}\}$\\ \hline
$8_1^3$&$4,2,2$&$\{\{2,2,2\},\{8,10,2,12,6,4\}\}$&$\{\{1,-3,2,-6\},\{3,-5,4,-1\},\{5,-2,6,-4\}\}$\\ \hline
$8_2^3$&$3\,1,2,2$&$\{\{2,2,2\},\{8,10,2,12,6,4\}\}$&$\{\{1,-3,2,-6\},\{3,-5,4,-1\},\{5,-2,6,-4\}\}$\\ \hline
$L10a142$&$5\,1,2,2$&$\{\{2,2,2\},\{8,10,2,12,6,4\}\}$&$\{\{1,-3,2,-6\},\{3,-5,4,-1\},\{5,-2,6,-4\}\}$\\ \hline
$L10a143$&$3\,2\,1,2,2$&$\{\{2,2,2\},\{8,10,2,12,6,4\}\}$&$\{\{1,-3,2,-6\},\{3,-5,4,-1\},\{5,-2,6,-4\}\}$\\ \hline
$L10a144$&$3\,3,2,2$&$\{\{2,2,2\},\{8,10,2,12,6,4\}\}$&$\{\{1,-3,2,-6\},\{3,-5,4,-1\},\{5,-2,6,-4\}\}$\\ \hline
$L10a145$&$6,2,2$&$\{\{2,2,2\},\{8,10,2,12,6,4\}\}$&$\{\{1,-3,2,-6\},\{3,-5,4,-1\},\{5,-2,6,-4\}\}$\\ \hline
$L10a146$&$(3,3)\,(2,2)$&$\{\{2,2,2\},\{8,10,2,12,6,4\}\}$&$\{\{1,-3,2,-6\},\{3,-5,4,-1\},\{5,-2,6,-4\}\}$\\ \hline
$L10a147$&$.3\,1:2\,0$&$\{\{2,2,2\},\{8,10,2,12,6,4\}\}$&$\{\{1,-3,2,-6\},\{3,-5,4,-1\},\{5,-2,6,-4\}\}$\\ \hline
$L10a148$&$.4:2\,0$&$\{\{2,2,2\},\{8,10,2,12,6,4\}\}$&$\{\{1,-3,2,-6\},\{3,-5,4,-1\},\{5,-2,6,-4\}\}$\\ \hline
$L10a159$&$3\,1,3\,1,2$&$\{\{2,2,2\},\{8,10,2,12,6,4\}\}$&$\{\{1,-3,2,-6\},\{3,-5,4,-1\},\{5,-2,6,-4\}\}$\\ \hline
$L10a160$&$3\,1,4,2$&$\{\{2,2,2\},\{8,10,2,12,6,4\}\}$&$\{\{1,-3,2,-6\},\{3,-5,4,-1\},\{5,-2,6,-4\}\}$\\ \hline
$L10a161$&$4,4,2$&$\{\{2,2,2\},\{8,10,2,12,6,4\}\}$&$\{\{1,-3,2,-6\},\{3,-5,4,-1\},\{5,-2,6,-4\}\}$\\ \hline
$L10a164$&$2\,0.2.2\,0.2\,0$&$\{\{2,2,2\},\{8,10,2,12,6,4\}\}$&$\{\{1,-3,2,-6\},\{3,-5,4,-1\},\{5,-2,6,-4\}\}$\\ \hline
&$101^*2:2$&$\{\{2,2,2\},\{8,10,2,12,6,4\}\}$&$\{\{1,-3,2,-6\},\{3,-5,4,-1\},\{5,-2,6,-4\}\}$\\ \hline
&$123^*$&$\{\{2,2,2\},\{8,10,2,12,6,4\}\}$&$\{\{1,-3,2,-6\},\{3,-5,4,-1\},\{5,-2,6,-4\}\}$\\ \hline
&$(2,2)\,1\,(3,3+)$&$\{\{2,2,2\},\{8,10,2,12,6,4\}\}$&$\{\{1,-3,2,-6\},\{3,-5,4,-1\},\{5,-2,6,-4\}\}$\\ \hline
&$(2,2)\,(3,3++)$&$\{\{2,2,2\},\{8,10,2,12,6,4\}\}$&$\{\{1,-3,2,-6\},\{3,-5,4,-1\},\{5,-2,6,-4\}\}$\\ \hline
&$(2,3\,1)\,(3,3)$&$\{\{2,2,2\},\{8,10,2,12,6,4\}\}$&$\{\{1,-3,2,-6\},\{3,-5,4,-1\},\{5,-2,6,-4\}\}$\\ \hline
&$(2,4)\,(3,3)$&$\{\{2,2,2\},\{8,10,2,12,6,4\}\}$&$\{\{1,-3,2,-6\},\{3,-5,4,-1\},\{5,-2,6,-4\}\}$\\ \hline
&$3\,1,3\,1,3\,1$&$\{\{2,2,2\},\{8,10,2,12,6,4\}\}$&$\{\{1,-3,2,-6\},\{3,-5,4,-1\},\{5,-2,6,-4\}\}$\\ \hline
&$3\,1,3\,1,4$&$\{\{2,2,2\},\{8,10,2,12,6,4\}\}$&$\{\{1,-3,2,-6\},\{3,-5,4,-1\},\{5,-2,6,-4\}\}$\\ \hline
&$3\,1,4,4$&$\{\{2,2,2\},\{8,10,2,12,6,4\}\}$&$\{\{1,-3,2,-6\},\{3,-5,4,-1\},\{5,-2,6,-4\}\}$\\ \hline
&$3\,1,6,2$&$\{\{2,2,2\},\{8,10,2,12,6,4\}\}$&$\{\{1,-3,2,-6\},\{3,-5,4,-1\},\{5,-2,6,-4\}\}$\\ \hline
&$3\,2\,1,3\,1,2$&$\{\{2,2,2\},\{8,10,2,12,6,4\}\}$&$\{\{1,-3,2,-6\},\{3,-5,4,-1\},\{5,-2,6,-4\}\}$\\ \hline
&$3\,2\,1,4,2$&$\{\{2,2,2\},\{8,10,2,12,6,4\}\}$&$\{\{1,-3,2,-6\},\{3,-5,4,-1\},\{5,-2,6,-4\}\}$\\ \hline
&$3\,2\,2\,1,2,2$&$\{\{2,2,2\},\{8,10,2,12,6,4\}\}$&$\{\{1,-3,2,-6\},\{3,-5,4,-1\},\{5,-2,6,-4\}\}$\\ \hline
&$3\,2\,3,2,2$&$\{\{2,2,2\},\{8,10,2,12,6,4\}\}$&$\{\{1,-3,2,-6\},\{3,-5,4,-1\},\{5,-2,6,-4\}\}$\\ \hline
&$(3,3\,2)\,(2,2)$&$\{\{2,2,2\},\{8,10,2,12,6,4\}\}$&$\{\{1,-3,2,-6\},\{3,-5,4,-1\},\{5,-2,6,-4\}\}$\\ \hline
&$3\,3,3\,1,2$&$\{\{2,2,2\},\{8,10,2,12,6,4\}\}$&$\{\{1,-3,2,-6\},\{3,-5,4,-1\},\{5,-2,6,-4\}\}$\\ \hline
&$3\,3,4,2$&$\{\{2,2,2\},\{8,10,2,12,6,4\}\}$&$\{\{1,-3,2,-6\},\{3,-5,4,-1\},\{5,-2,6,-4\}\}$\\ \hline
&$3\,4\,1,2,2$&$\{\{2,2,2\},\{8,10,2,12,6,4\}\}$&$\{\{1,-3,2,-6\},\{3,-5,4,-1\},\{5,-2,6,-4\}\}$\\ \hline
&$(3,5)\,(2,2)$&$\{\{2,2,2\},\{8,10,2,12,6,4\}\}$&$\{\{1,-3,2,-6\},\{3,-5,4,-1\},\{5,-2,6,-4\}\}$\\ \hline
&$3\,5,2,2$&$\{\{2,2,2\},\{8,10,2,12,6,4\}\}$&$\{\{1,-3,2,-6\},\{3,-5,4,-1\},\{5,-2,6,-4\}\}$\\ \hline
&$4,4,4$&$\{\{2,2,2\},\{8,10,2,12,6,4\}\}$&$\{\{1,-3,2,-6\},\{3,-5,4,-1\},\{5,-2,6,-4\}\}$\\ \hline
&$5\,1,3\,1,2$&$\{\{2,2,2\},\{8,10,2,12,6,4\}\}$&$\{\{1,-3,2,-6\},\{3,-5,4,-1\},\{5,-2,6,-4\}\}$\\ \hline
&$5\,1,4,2$&$\{\{2,2,2\},\{8,10,2,12,6,4\}\}$&$\{\{1,-3,2,-6\},\{3,-5,4,-1\},\{5,-2,6,-4\}\}$\\ \hline
&$5\,2\,1,2,2$&$\{\{2,2,2\},\{8,10,2,12,6,4\}\}$&$\{\{1,-3,2,-6\},\{3,-5,4,-1\},\{5,-2,6,-4\}\}$\\ \hline
&$5\,3,2,2$&$\{\{2,2,2\},\{8,10,2,12,6,4\}\}$&$\{\{1,-3,2,-6\},\{3,-5,4,-1\},\{5,-2,6,-4\}\}$\\ \hline
&$6^*2.2\,0.2.2\,0.2.2\,0$&$\{\{2,2,2\},\{8,10,2,12,6,4\}\}$&$\{\{1,-3,2,-6\},\{3,-5,4,-1\},\{5,-2,6,-4\}\}$\\ \hline
&$6^*2.2:(2,2).2\,0$&$\{\{2,2,2\},\{8,10,2,12,6,4\}\}$&$\{\{1,-3,2,-6\},\{3,-5,4,-1\},\{5,-2,6,-4\}\}$\\ \hline
&$6^*2.2.2.3\,0.2$&$\{\{2,2,2\},\{8,10,2,12,6,4\}\}$&$\{\{1,-3,2,-6\},\{3,-5,4,-1\},\{5,-2,6,-4\}\}$\\ \hline
&$6^*2.2.2:3\,1$&$\{\{2,2,2\},\{8,10,2,12,6,4\}\}$&$\{\{1,-3,2,-6\},\{3,-5,4,-1\},\{5,-2,6,-4\}\}$\\ \hline
&$6^*2.2.2:4$&$\{\{2,2,2\},\{8,10,2,12,6,4\}\}$&$\{\{1,-3,2,-6\},\{3,-5,4,-1\},\{5,-2,6,-4\}\}$\\ \hline
&$6^*2.3\,1.2:2$&$\{\{2,2,2\},\{8,10,2,12,6,4\}\}$&$\{\{1,-3,2,-6\},\{3,-5,4,-1\},\{5,-2,6,-4\}\}$\\ \hline
&$6^*2:.3\,2\,1$&$\{\{2,2,2\},\{8,10,2,12,6,4\}\}$&$\{\{1,-3,2,-6\},\{3,-5,4,-1\},\{5,-2,6,-4\}\}$\\ \hline
&$6^*2.3.2.2.2$&$\{\{2,2,2\},\{8,10,2,12,6,4\}\}$&$\{\{1,-3,2,-6\},\{3,-5,4,-1\},\{5,-2,6,-4\}\}$\\ \hline
&$6^*2:.(3,3)$&$\{\{2,2,2\},\{8,10,2,12,6,4\}\}$&$\{\{1,-3,2,-6\},\{3,-5,4,-1\},\{5,-2,6,-4\}\}$\\ \hline
&$6^*2:.3\,3$&$\{\{2,2,2\},\{8,10,2,12,6,4\}\}$&$\{\{1,-3,2,-6\},\{3,-5,4,-1\},\{5,-2,6,-4\}\}$\\ \hline
&$6^*2.3:3\,1$&$\{\{2,2,2\},\{8,10,2,12,6,4\}\}$&$\{\{1,-3,2,-6\},\{3,-5,4,-1\},\{5,-2,6,-4\}\}$\\ \hline
&$6^*2.3:4$&$\{\{2,2,2\},\{8,10,2,12,6,4\}\}$&$\{\{1,-3,2,-6\},\{3,-5,4,-1\},\{5,-2,6,-4\}\}$\\ \hline
&$6^*2.4.2:2$&$\{\{2,2,2\},\{8,10,2,12,6,4\}\}$&$\{\{1,-3,2,-6\},\{3,-5,4,-1\},\{5,-2,6,-4\}\}$\\ \hline
&$6^*2:.5\,1$&$\{\{2,2,2\},\{8,10,2,12,6,4\}\}$&$\{\{1,-3,2,-6\},\{3,-5,4,-1\},\{5,-2,6,-4\}\}$\\ \hline
&$6^*2:.6$&$\{\{2,2,2\},\{8,10,2,12,6,4\}\}$&$\{\{1,-3,2,-6\},\{3,-5,4,-1\},\{5,-2,6,-4\}\}$\\ \hline
&$6^*3\,1.2.2:2$&$\{\{2,2,2\},\{8,10,2,12,6,4\}\}$&$\{\{1,-3,2,-6\},\{3,-5,4,-1\},\{5,-2,6,-4\}\}$\\ \hline
&$6^*3\,1:.3\,1$&$\{\{2,2,2\},\{8,10,2,12,6,4\}\}$&$\{\{1,-3,2,-6\},\{3,-5,4,-1\},\{5,-2,6,-4\}\}$\\ \hline
&$6^*3\,1.3:2$&$\{\{2,2,2\},\{8,10,2,12,6,4\}\}$&$\{\{1,-3,2,-6\},\{3,-5,4,-1\},\{5,-2,6,-4\}\}$\\ \hline
&$6^*3\,1:.4$&$\{\{2,2,2\},\{8,10,2,12,6,4\}\}$&$\{\{1,-3,2,-6\},\{3,-5,4,-1\},\{5,-2,6,-4\}\}$\\ \hline
&$6^*3.3\,0.3$&$\{\{2,2,2\},\{8,10,2,12,6,4\}\}$&$\{\{1,-3,2,-6\},\{3,-5,4,-1\},\{5,-2,6,-4\}\}$\\ \hline
&$6^*3.3\,0::3\,0$&$\{\{2,2,2\},\{8,10,2,12,6,4\}\}$&$\{\{1,-3,2,-6\},\{3,-5,4,-1\},\{5,-2,6,-4\}\}$\\ \hline
&$6^*3.3.3$&$\{\{2,2,2\},\{8,10,2,12,6,4\}\}$&$\{\{1,-3,2,-6\},\{3,-5,4,-1\},\{5,-2,6,-4\}\}$\\ \hline
&$6^*3.3.3\,0$&$\{\{2,2,2\},\{8,10,2,12,6,4\}\}$&$\{\{1,-3,2,-6\},\{3,-5,4,-1\},\{5,-2,6,-4\}\}$\\ \hline
&$6,4,2$&$\{\{2,2,2\},\{8,10,2,12,6,4\}\}$&$\{\{1,-3,2,-6\},\{3,-5,4,-1\},\{5,-2,6,-4\}\}$\\ \hline
&$6^*4.2.2:2$&$\{\{2,2,2\},\{8,10,2,12,6,4\}\}$&$\{\{1,-3,2,-6\},\{3,-5,4,-1\},\{5,-2,6,-4\}\}$\\ \hline
&$6^*4.3:2$&$\{\{2,2,2\},\{8,10,2,12,6,4\}\}$&$\{\{1,-3,2,-6\},\{3,-5,4,-1\},\{5,-2,6,-4\}\}$\\ \hline
&$7\,1,2,2$&$\{\{2,2,2\},\{8,10,2,12,6,4\}\}$&$\{\{1,-3,2,-6\},\{3,-5,4,-1\},\{5,-2,6,-4\}\}$\\ \hline
&$8,2,2$&$\{\{2,2,2\},\{8,10,2,12,6,4\}\}$&$\{\{1,-3,2,-6\},\{3,-5,4,-1\},\{5,-2,6,-4\}\}$\\ \hline
&$8^*2:2\,0.2\,0:.2\,0$&$\{\{2,2,2\},\{8,10,2,12,6,4\}\}$&$\{\{1,-3,2,-6\},\{3,-5,4,-1\},\{5,-2,6,-4\}\}$\\ \hline
&$8^*2.2\,0:.2\,0:.2\,0$&$\{\{2,2,2\},\{8,10,2,12,6,4\}\}$&$\{\{1,-3,2,-6\},\{3,-5,4,-1\},\{5,-2,6,-4\}\}$\\ \hline
&$8^*2:2\,0:2:2\,0$&$\{\{2,2,2\},\{8,10,2,12,6,4\}\}$&$\{\{1,-3,2,-6\},\{3,-5,4,-1\},\{5,-2,6,-4\}\}$\\ \hline
&$8^*2:2:2:2$&$\{\{2,2,2\},\{8,10,2,12,6,4\}\}$&$\{\{1,-3,2,-6\},\{3,-5,4,-1\},\{5,-2,6,-4\}\}$\\ \hline
&$8^*2.2.2:.2$&$\{\{2,2,2\},\{8,10,2,12,6,4\}\}$&$\{\{1,-3,2,-6\},\{3,-5,4,-1\},\{5,-2,6,-4\}\}$\\ \hline
&$8^*2.2:3$&$\{\{2,2,2\},\{8,10,2,12,6,4\}\}$&$\{\{1,-3,2,-6\},\{3,-5,4,-1\},\{5,-2,6,-4\}\}$\\ \hline
&$8^*2.2:3\,0$&$\{\{2,2,2\},\{8,10,2,12,6,4\}\}$&$\{\{1,-3,2,-6\},\{3,-5,4,-1\},\{5,-2,6,-4\}\}$\\ \hline
&$9^*2\,0.2\,0:::2\,0$&$\{\{2,2,2\},\{8,10,2,12,6,4\}\}$&$\{\{1,-3,2,-6\},\{3,-5,4,-1\},\{5,-2,6,-4\}\}$\\ \hline
&$9^*2:.2\,0::.2\,0$&$\{\{2,2,2\},\{8,10,2,12,6,4\}\}$&$\{\{1,-3,2,-6\},\{3,-5,4,-1\},\{5,-2,6,-4\}\}$\\ \hline
&$9^*2::2::2$&$\{\{2,2,2\},\{8,10,2,12,6,4\}\}$&$\{\{1,-3,2,-6\},\{3,-5,4,-1\},\{5,-2,6,-4\}\}$\\ \hline
\end{tabular}

\normalsize

\bigskip

The next question is: which multi-component links can be represented by meander diagrams. Same as before, all multi-component meander links belong to their corresponding families: if $p$ is a twist in some meander link, its replacement by $p+2$ results in the meander link with the same number of components.

\section{Sum of meander knots and links}

\begin{figure}[th]
\centerline{\psfig{file=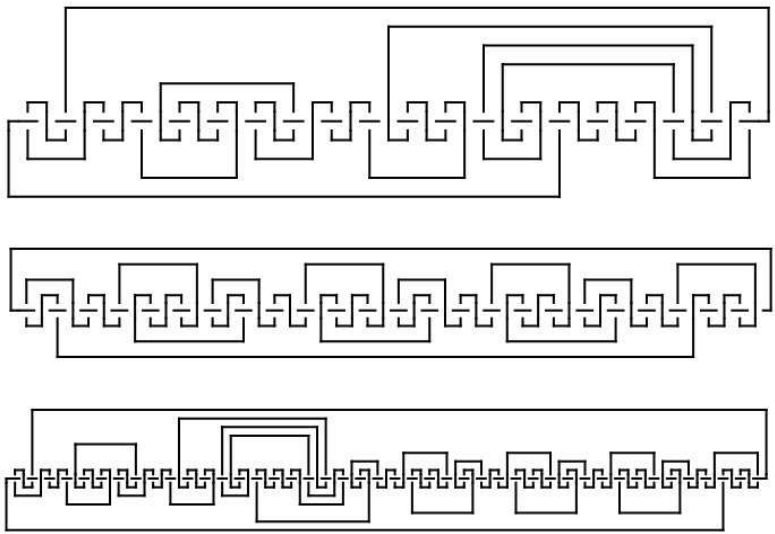,width=4.2in}} \vspace*{8pt}
\caption{Sum of a 39-crossing meander knot and 48-crossing meander link giving the 87-crossing meander knot.
\label{sum}}
\end{figure}

\begin{figure}[th]
\centerline{\psfig{file=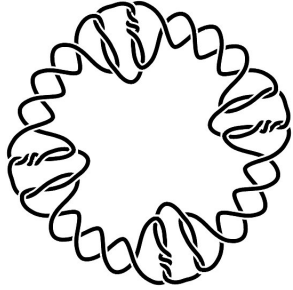,width=1.5in}} \vspace*{8pt}
\caption{A 48-crossing meander link obtained by iterative sum.
\label{sum1}}
\end{figure}

For two open meander sequences we can define a {\it sum or concatenation}: the operation of joining their Dyck words and connecting the second loose end of the first with the first loose end of the second and making a closure in order to obtain a meander knot or link ($KL$) diagram (Fig. \ref{sum}). The same definition extends to meander knots and links where we concatenate the meander parts of their Gauss codes. From the parity reasons, the sum of two meander knot diagrams or the sum of two meander link diagrams is a meander link diagram, and the sum of a meander knot diagram and meander link diagram or {\it vice versa} is a meander knot diagram. Sum of a meander knot diagram and its mirror image is a 2-component unlink.  An interesting question is to analyze sums of meander $KL$ diagrams from the point of view of tangle operations which can be recognized from Conway symbols of $KL$s and their sums. Some artistic meander knots and links can be obtained by iterative sum (Fig. \ref{sum1}).

\section{Semi-meander or ordered Gauss code knots}

\begin{figure}[th]
\centerline{\psfig{file=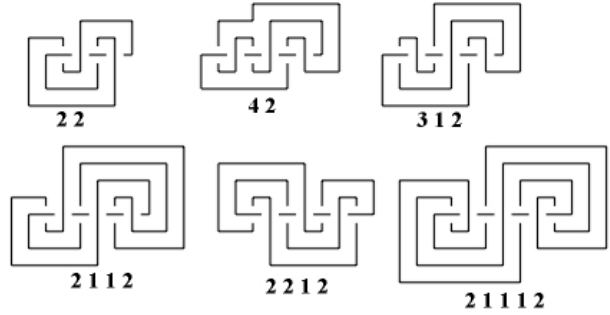,width=3.2in}} \vspace*{8pt}
\caption{Semi-meander knots with $n\le 7$ crossings which are not meander knots.
\label{semi01}}
\end{figure}

\begin{figure}[th]
\centerline{\psfig{file=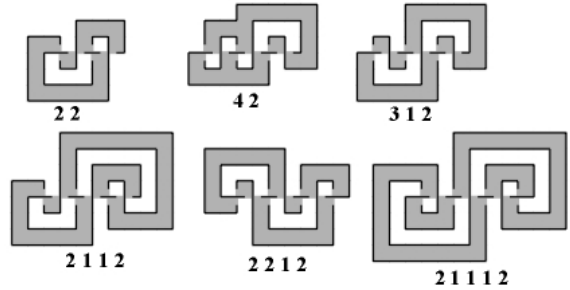,width=3.0in}} \vspace*{8pt}
\caption{Checker-board colorings of knots from Fig. \ref{semi01}.
\label{semi02}}
\end{figure}

In the case of meanders the axis of a meander is infinite. If the axis is finite, as the result we obtain semi-meanders, where after the end(s) of the axis meander curve can pass from one side of the axis to the other without crossing the axis. A Gauss code of every knot depends from the choice of the initial (basic) point belonging to some arc and from the orientation of the knot. This means that every rotation or reversal of a sequence of length $2n$ representing a Gauss code of a knot with $n$ crossings represents the same (non-oriented) knot. A Gauss code will be called {\it ordered} if the absolute value of the first part of its Gauss code is the sequence $1,2,\ldots ,n$. An alternating knot will be called {\it Gauss code ordered} (OGC) or {\it semi-meander} knot if it has at least one minimal diagram with ordered Gauss code. The name {\it semi-meander knot} follows from the fact that the shadow of such a knot represents a meander or semi-meander. It is clear that every meander knot is GCO, and that meander knots represent the proper subset of GCO knots. For OGC knots there is no parity restriction to the number of crossings, so there exist OGC knots which are not meander knots: OGC knots with an even number of crossings. Moreover, some OGC knots with an odd number of crossings are not meander knots, e.g., knot $7_6=2\,2\,1\,2$ which has two minimal diagrams, and among them only one is OGC diagram with ordered Gauss code $\{1,-2,3,-4,5,-6,7,-5,4,-1,2,-7,6,-3\}$. Every OGC diagram is completely determined by the second half of its ordered Gauss code, which will be called short Gauss code. Fig. \ref{semi01} shows all semi-meander knots with $n\le 7$ crossings which are not meander knots, and Fig. \ref{semi02} their checker-board colorings. All alternating OGC knots are members of their corresponding knot families.

First alternating knots which are not OGC knots occur for $n=8$, and they are knots $8_{16}=.2.2\,0$ and $8_{18}=8^*$, for $n=9$ there are 6 knots that are not OGC: $9_{29}=.2.2\,0.2$, $9_{32}=.2\,1.2\,0$, $9_{33}=.2\,1.2$, $9_{34}=8^*2\,0$, $9_{40}=9^*$, $9_41=2\,0:2\,0:2\,0$, {\it etc.} In the following table is given the number of alternating knots with $n$ crossings ($3\le n\le 12$), the number of OGC knots, and the number of OGC knots that are not meander knots.

In the following table are given all alternating OGC knots with $4\le n\le 9$ crossings which are not meander knots. Every knot is given by its standard symbol, Conway symbol, OGC diagram Dowker-Thistlethwaite code and short OGC Gauss code.

\bigskip

\tiny

\noindent \begin{tabular}{|c|c|c|c|}   \hline
$4_1$&$2\,2$&$\{\{4\},\{4,6,8,2\}\}$&$\{-2,1,-4,3\}$\\ \hline
$6_1$&$4\,2$&$\{\{6\},\{4,8,12,10,2,6\}\}$&$\{-4,3,-2,1,-6,5\}$\\ \hline
$6_2$&$3\,1\,2$&$\{\{6\},\{4,8,10,12,2,6\}\}$&$\{-4,1,-2,3,-6,5\}$\\ \hline
$6_3$&$2\,1\,1\,2$&$\{\{6\},\{4,8,10,2,12,6\}\}$&$\{4,-1,2,-5,6,-3\}$\\ \hline
$7_6$&$2\,2\,1\,2$&$\{\{7\},\{4,8,12,2,14,6,10\}\}$&$\{-5,4,-1,2,-7,6,-3\}$\\ \hline
$7_7$&$2\,1\,1\,1\,2$&$\{\{7\},\{4,8,12,14,2,6,10\}\}$&$\{-5,2,-1,4,-7,6,-3\}$\\ \hline
$8_1$&$6\,2$&$\{\{8\},\{4,10,16,14,12,2,8,6\}\}$&$\{-6,5,-4,3,-2,1,-8,7\}$\\ \hline
$8_2$&$5\,1\,2$&$\{\{8\},\{4,10,12,14,16,2,6,8\}\}$&$\{-6,1,-2,3,-4,5,-8,7\}$\\ \hline
$8_3$&$4\,4$&$\{\{8\},\{6,12,10,16,14,4,2,8\}\}$&$\{-4,3,-2,1,-8,7,-6,5\}$\\ \hline
$8_4$&$4\,1\,3$&$\{\{8\},\{6,10,12,16,14,4,2,8\}\}$&$\{-4,3,-2,1,-8,5,-6,7\}$\\ \hline
$8_5$&$3,3,2$&$\{\{8\},\{6,8,12,2,14,16,4,10\}\}$&$\{-4,5,-6,1,-2,3,-8,7\}$\\ \hline
$8_6$&$3\,3\,2$&$\{\{8\},\{4,10,16,12,14,2,8,6\}\}$&$\{-6,5,-4,1,-2,3,-8,7\}$\\ \hline
$8_7$&$4\,1\,1\,2$&$\{\{8\},\{4,10,12,14,2,16,6,8\}\}$&$\{4,-1,2,-5,6,-7,8,-3\}$\\ \hline
$8_8$&$2\,3\,1\,2$&$\{\{8\},\{4,8,12,2,16,14,6,10\}\}$&$\{6,-5,4,-1,2,-7,8,-3\}$\\ \hline
$8_9$&$3\,1\,1\,3$&$\{\{8\},\{6,10,12,14,16,4,2,8\}\}$&$\{-4,1,-2,3,-8,5,-6,7\}$\\ \hline
$8_{10}$&$2\,1,3,2$&$\{\{8\},\{4,8,12,2,14,16,6,10\}\}$&$\{6,-1,2,-7,8,-3,4,-5\}$\\ \hline
$8_{11}$&$3\,2\,1\,2$&$\{\{8\},\{4,10,12,14,16,2,8,6\}\}$&$\{-6,1,-4,3,-2,5,-8,7\}$\\ \hline
$8_{12}$&$2\,2\,2\,2$&$\{\{8\},\{4,8,16,12,2,14,6,10\}\}$&$\{-6,5,-2,1,-8,7,-4,3\}$\\ \hline
$8_{13}$&$3\,1\,1\,1\,2$&$\{\{8\},\{4,10,12,14,2,16,8,6\}\}$&$\{4,-1,2,-7,6,-5,8,-3\}$\\ \hline
$8_{14}$&$2\,2\,1\,1\,2$&$\{\{8\},\{4,8,12,16,2,14,6,10\}\}$&$\{-6,5,-2,1,-4,7,-8,3\}$\\ \hline
$8_{15}$&$2\,1,2\,1,2$&$\{\{8\},\{4,8,14,2,12,16,6,10\}\}$&$\{4,-5,8,-1,2,-7,6,-3\}$\\ \hline
$8_{17}$&$.2.2$&$\{\{8\},\{6,8,12,14,4,16,2,10\}\}$&$\{2,-5,6,-1,8,-3,4,-7\}$\\ \hline
\end{tabular}

\noindent \begin{tabular}{|c|c|c|c|}   \hline
$9_8$&$2\,4\,1\,2$&$\{\{9\},\{4,8,14,2,18,16,6,12,10\}\}$&$\{-7,6,-5,4,-1,2,-9,8,-3\}$\\ \hline
$9_{11}$&$4\,1\,2\,2$&$\{\{9\},\{4,10,18,14,16,2,6,8,12\}\}$&$\{5,-2,1,-6,7,-8,9,-4,3\}$\\ \hline
$9_{12}$&$4\,2\,1\,2$&$\{\{9\},\{4,10,16,14,2,18,8,6,12\}\}$&$\{-5,4,-1,2,-9,8,-7,6,-3\}$\\ \hline
$9_{14}$&$4\,1\,1\,1\,2$&$\{\{9\},\{4,10,16,14,18,2,8,6,12\}\}$&$\{-7,4,-3,2,-1,6,-9,8,-5\}$\\ \hline
$9_{15}$&$2\,3\,2\,2$&$\{\{9\},\{4,8,18,14,2,16,6,12,10\}\}$&$\{7,-6,5,-2,1,-8,9,-4,3\}$\\ \hline
$9_{17}$&$2\,1\,3\,1\,2$&$\{\{9\},\{4,10,14,16,18,2,6,12,8\}\}$&$\{-7,2,-1,4,-5,6,-9,8,-3\}$\\ \hline
$9_{19}$&$2\,3\,1\,1\,2$&$\{\{9\},\{4,8,14,18,2,16,6,12,10\}\}$&$\{7,-6,5,-2,1,-4,9,-8,3\}$\\ \hline
$9_{20}$&$3\,1\,2\,1\,2$&$\{\{9\},\{4,10,14,16,2,18,8,6,12\}\}$&$\{-5,4,-1,2,-9,6,-7,8,-3\}$\\ \hline
$9_{21}$&$3\,1\,1\,2\,2$&$\{\{9\},\{4,10,18,16,14,2,6,8,12\}\}$&$\{5,-2,1,-8,7,-6,9,-4,3\}$\\ \hline
$9_{22}$&$2\,1\,1,3,2$&$\{\{9\},\{4,8,14,18,2,16,6,10,12\}\}$&$\{-5,6,-7,2,-1,4,-9,8,-3\}$\\ \hline
$9_{24}$&$2\,1,3,2+$&$\{\{9\},\{4,8,14,2,18,16,6,10,12\}\}$&$\{7,-2,1,-8,9,-6,3,-4,5\}$\\ \hline
$9_{25}$&$2\,2,2\,1,2$&$\{\{9\},\{4,8,16,2,14,18,10,6,12\}\}$&$\{-7,6,-3,2,-9,8,-1,4,-5\}$\\ \hline
$9_{26}$&$3\,1\,1\,1\,1\,2$&$\{\{9\},\{4,10,14,16,18,2,8,6,12\}\}$&$\{7,-2,3,-4,1,-6,9,-8,5\}$\\ \hline
$9_{27}$&$2\,1\,2\,1\,1\,2$&$\{\{9\},\{4,10,14,16,2,18,6,12,8\}\}$&$\{7,-4,1,-2,5,-6,9,-8,3\}$\\ \hline
$9_{28}$&$2\,1,2\,1,2+$&$\{\{9\},\{4,8,14,2,18,16,10,6,12\}\}$&$\{7,-2,3,-8,9,-6,1,-4,5\}$\\ \hline
$9_{30}$&$2\,1\,1,2\,1,2$&$\{\{9\},\{4,8,14,2,16,18,10,6,12\}\}$&$\{7,-2,3,-6,9,-8,1,-4,5\}$\\ \hline
$9_{31}$&$2\,1\,1\,1\,1\,1\,2$&$\{\{9\},\{4,10,16,14,2,18,6,12,8\}\}$&$\{-7,4,-1,2,-5,8,-9,6,-3\}$\\ \hline
$9_{36}$&$2\,2,3,2$&$\{\{9\},\{4,8,18,14,2,16,6,10,12\}\}$&$\{-7,6,-1,2,-9,8,-3,4,-5\}$\\ \hline
$9_{37}$&$2\,1,2\,1,3$&$\{\{9\},\{4,10,14,18,16,2,6,12,8\}\}$&$\{-3,2,-9,6,-5,4,-1,8,-7\}$\\ \hline
$9_{39}$&$2:2:2\,0$&$\{\{9\},\{6,10,14,18,16,2,8,4,12\}\}$&$\{5,-2,1,-6,9,-4,3,-8,7\}$\\ \hline
\end{tabular}

\normalsize

\bigskip

In the following table is given the number of non-alternating knots for a given number of crossings $n$ ($8\le n\le 12$), the number of OGC knots, and the number of OGC knots which are not meander knots:

\bigskip

\small

\begin{tabular}{|c|c|c|c|}   \hline
$n$ & No. of knots & OGC & OGC non-meander  \\  \hline
8 & 3 & 3 & 3 \\ \hline
9 & 8 & 7 & 6 \\ \hline
10 & 42 & 39 & 39 \\ \hline
11 & 187 & 141 & 73 \\ \hline
12 & 888 & 636 & 636 \\ \hline

\end{tabular}

\normalsize

\bigskip

In the following table are given OGC diagrams of all non-alternating knots with $n\le 10$ crossings which don't have a meander diagram with the number of crossings equal to the crossing number. Every knot is given by its standard symbol, Conway symbol, OGC diagram Dowker-Thistlethwaite code and short OGC Gauss code.

\bigskip

\tiny

\noindent \begin{tabular}{|c|c|c|c|}   \hline
$8_{19}$&$3,3,-2$&$\{\{8\},\{6,8,-12,2,14,16,-4,10\}\}$&$\{-4,5,-6,1,-2,3,8,-7\}\}$\\ \hline
$8_{20}$&$3,2\,1,-2$&$\{\{8\},\{4,8,-12,2,14,16,-6,10\}\}$&$\{2,-3,-6,7,-8,1,-4,5\}\}$\\ \hline
$8_{21}$&$2\,1,2\,1,-2$&$\{\{8\},\{4,8,-14,2,12,16,-6,10\}\}$&$\{4,-5,8,-1,2,7,-6,-3\}\}$\\ \hline
$9_{42}$&$2\,2,3,-2$&$\{\{9\},\{-4,-8,-18,14,-2,-16,6,-10,-12\}\}$&$\{-7,6,1,-2,-9,8,-3,4,-5\}$\\ \hline
$9_{43}$&$2\,1\,1,3,-2$&$\{\{9\},\{4,8,-14,18,2,16,-6,10,12\}\}$&$\{-5,6,-7,2,-1,4,9,-8,-3\}$\\ \hline
$9_{44}$&$2\,2,2\,1,-2$&$\{\{9\},\{4,8,-16,2,14,18,10,-6,12\}\}$&$\{-7,6,3,-2,-9,8,-1,4,-5\}$\\ \hline
$9_{45}$&$2\,1\,1,2\,1,-2$&$\{\{9\},\{4,8,-14,2,16,18,10,-6,12\}\}$&$\{7,2,-3,-6,9,-8,1,-4,5\}$\\ \hline
$9_{46}$&$3,3,-3$&$\{\{9\},\{8,-12,16,14,18,-4,-2,6,10\}\}$&$\{3,-2,1,-6,5,-4,-9,8,-7\}$\\ \hline
$9_{48}$&$2\,1,2\,1,-3$&$\{\{9\},\{4,10,14,-18,-16,2,6,12,-8\}\}$&$\{-3,2,-9,-6,5,-4,-1,8,-7\}$\\ \hline
$9_{49}$&$-2\,0:-2\,0:-2\,0$&$\{\{9\},\{-6,10,14,-18,16,2,8,4,12\}\}$&$\{5,-2,1,-6,9,4,-3,-8,7\}$\\ \hline
$10_{124}$&$5,3,-2$&$\{\{10\},\{6,8,-14,2,16,18,20,-4,10,12\}\}$&$\{-6,7,-8,1,-2,3,-4,5,10,-9\}$\\ \hline
$10_{125}$&$5,2\,1,-2$&$\{\{10\},\{4,8,-14,2,16,18,20,-6,10,12\}\}$&$\{2,-3,-6,7,-8,9,-10,1,-4,5\}$\\ \hline
$10_{126}$&$4\,1,3,-2$&$\{\{10\},\{-6,-8,14,-2,-16,-18,4,-20,-10,-12\}\}$&$\{6,1,-2,-7,8,-9,10,-3,4,-5\}$\\ \hline
$10_{127}$&$4\,1,2\,1,-2$&$\{\{10\},\{4,8,-16,2,14,18,20,-6,10,12\}\}$&$\{6,-7,10,-1,2,-3,4,9,-8,-5\}$\\ \hline
$10_{128}$&$3\,2,3,-2$&$\{\{10\},\{6,8,-14,2,16,18,20,-4,12,10\}\}$&$\{-6,7,-8,1,-4,3,-2,5,10,-9\}$\\ \hline
$10_{129}$&$3\,2,2\,1,-2$&$\{\{10\},\{4,8,-14,2,18,16,20,-6,10,12\}\}$&$\{2,-3,-8,7,-6,9,-10,1,-4,5\}$\\ \hline
$10_{130}$&$3\,1\,1,3,-2$&$\{\{10\},\{-6,-8,14,-2,-16,-18,4,-20,-12,-10\}\}$&$\{6,1,-2,-9,8,-7,10,-3,4,-5\}$\\ \hline
$10_{131}$&$3\,1\,1,2\,1,-2$&$\{\{10\},\{4,8,-16,2,14,20,18,-6,10,12\}\}$&$\{6,-7,10,-3,2,-1,4,9,-8,-5\}$\\ \hline
$10_{132}$&$2\,3,3,-2$&$\{\{10\},\{-4,-10,-20,14,-2,-16,-18,8,-12,-6\}\}$&$\{8,-7,6,1,-2,-9,10,-3,4,-5\}$\\ \hline
$10_{133}$&$2\,3,2\,1,-2$&$\{\{10\},\{4,8,-18,2,16,14,20,10,-6,12\}\}$&$\{8,-7,6,3,-2,-9,10,-1,4,-5\}$\\ \hline
$10_{134}$&$2\,2\,1,3,-2$&$\{\{10\},\{-4,-10,-20,14,-2,-16,-18,6,-12,-8\}\}$&$\{-8,7,2,-1,-6,9,-10,3,-4,5\}$\\ \hline
$10_{135}$&$2\,2\,1,2\,1,-2$&$\{\{10\},\{4,8,-16,2,18,14,20,10,-6,12\}\}$&$\{-8,7,2,-3,-6,9,-10,1,-4,5\}$\\ \hline
$10_{136}$&$2\,2,2\,2,-2$&$\{\{10\},\{4,8,20,-16,2,14,18,10,-6,12\}\}$&$\{-6,5,-10,9,-2,1,8,-7,-4,3\}$\\ \hline
$10_{137}$&$2\,2,2\,1\,1,-2$&$\{\{10\},\{4,8,-16,20,2,14,18,10,-6,12\}\}$&$\{-6,5,-10,7,-2,1,8,-9,-4,3\}$\\ \hline
$10_{138}$&$2\,1\,1,2\,1\,1,-2$&$\{\{10\},\{4,8,-16,20,2,14,18,-6,10,12\}\}$&$\{8,3,-2,-7,10,-9,4,-1,6,-5\}$\\ \hline
$10_{139}$&$4,3,-2\,-1$&$\{\{10\},\{4,10,-14,-16,2,-18,-20,-6,-8,-12\}\}$&$\{4,-5,6,1,-2,-7,8,-9,10,3\}$\\ \hline
$10_{140}$&$4,3,-3$&$\{\{10\},\{-8,-10,16,14,-2,18,20,6,4,12\}\}$&$\{-4,3,-2,1,-8,9,-10,-5,6,-7\}$\\ \hline
\end{tabular}

\noindent \begin{tabular}{|c|c|c|c|}   \hline
$10_{141}$&$4,2\,1,-3$&$\{\{10\},\{-4,-10,-14,-16,-2,18,20,-6,-8,12\}\}$&$\{-4,5,-10,-1,2,-3,-6,7,-8,9\}$\\ \hline
$10_{142}$&$3\,1,3,-2\,-1$&$\{\{10\},\{4,10,-14,-16,2,-18,-20,-8,-6,-12\}\}$&$\{4,-5,6,1,-2,-9,8,-7,10,3\}$\\ \hline
$10_{143}$&$3\,1,3,-3$&$\{\{10\},\{-8,-10,14,16,-2,18,20,6,4,12\}\}$&$\{-4,1,-2,3,-8,9,-10,-5,6,-7\}$\\ \hline
$10_{144}$&$3\,1,2\,1,-3$&$\{\{10\},\{-4,-10,-14,-16,-2,18,20,-8,-6,12\}\}$&$\{-4,5,-10,-1,2,-3,-8,7,-6,9\}$\\ \hline
$10_{145}$&$2\,2,3,-2\,-1$&$\{\{10\},\{-4,-10,14,20,18,-2,16,6,12,8\}\}$&$\{-8,-1,2,-7,6,-5,-10,9,4,-3\}$\\ \hline
$10_{146}$&$2\,2,2\,1,-3$&$\{\{10\},\{4,10,14,-20,-18,2,16,6,12,-8\}\}$&$\{8,-7,2,-1,-6,5,-4,-9,10,-3\}$\\ \hline
$10_{147}$&$2\,1\,1,3,-3$&$\{\{10\},\{4,12,16,14,-18,2,20,6,-10,-8\}\}$&$\{8,-3,2,-1,-6,5,-4,-9,10,-7\}$\\ \hline
$10_{148}$&$(3,2)\,(3,-2)$&$\{\{10\},\{6,8,12,2,-16,4,18,20,-10,14\}\}$&$\{6,-7,8,-1,2,9,-10,-3,4,-5\}$\\ \hline
$10_{149}$&$(3,2)\,(2\,1,-2)$&$\{\{10\},\{4,8,-18,2,14,16,20,10,-6,12\}\}$&$\{6,-7,10,-1,2,9,-8,-3,4,-5\}$\\ \hline
$10_{150}$&$(2\,1,2)\,(3,-2)$&$\{\{10\},\{4,8,18,2,14,16,-20,10,6,-12\}\}$&$\{6,-7,8,-3,2,9,-10,-1,4,-5\}$\\ \hline
$10_{151}$&$(2\,1,2)\,(2\,1,-2)$&$\{\{10\},\{4,8,-18,2,14,20,16,10,-6,12\}\}$&$\{6,-7,10,-3,2,9,-8,-1,4,-5\}$\\ \hline
$10_{152}$&$(3,2)\,-(3,2)$&$\{\{10\},\{6,8,12,2,-16,4,-18,-20,-10,-14\}\}$&$\{-6,7,-8,-1,2,9,-10,-3,4,-5\}$\\ \hline
$10_{153}$&$(3,2)\,-(2\,1,2)$&$\{\{10\},\{-4,-8,-18,-2,14,16,20,10,-6,12\}\}$&$\{-6,7,-10,-1,2,9,-8,-3,4,-5\}$\\ \hline
$10_{154}$&$(2\,1,2)\,-(2\,1,2)$&$\{\{10\},\{-4,-8,-18,-2,14,20,16,10,-6,12\}\}$&$\{-6,7,-10,-3,2,9,-8,-1,4,-5\}$\\ \hline
$10_{155}$&$-3:2:2$&$\{\{10\},\{6,10,14,16,18,4,-20,2,8,-12\}\}$&$\{6,-1,2,-3,8,-9,-4,5,10,-7\}$\\ \hline
$10_{156}$&$-3:2:2\,0$&$\{\{10\},\{6,10,14,-18,-16,4,20,2,-8,12\}\}$&$\{2,-5,6,-1,10,-3,4,9,-8,7\}$\\ \hline
$10_{158}$&$-3\,0:2:2$&$\{\{10\},\{-6,10,16,-20,14,2,18,4,8,12\}\}$&$\{6,-1,2,-5,10,-9,-4,3,8,-7\}$\\ \hline
$10_{159}$&$-3\,0:2:2\,0$&$\{\{10\},\{6,10,14,-16,-18,4,20,2,-8,12\}\}$&$\{2,-5,6,-1,10,-3,4,7,-8,9\}$\\ \hline
$10_{160}$&$-3\,0:2\,0:2\,0$&$\{\{10\},\{-6,-10,16,-20,14,-2,-18,4,8,-12\}\}$&$\{-6,1,-2,5,10,-9,4,-3,8,-7\}$\\ \hline
$10_{161}$&$3:-2\,0:-2\,0$&$\{\{10\},\{6,-10,16,20,14,-2,-18,4,8,-12\}\}$&$\{6,-1,2,-5,-10,9,4,-3,-8,7\}$\\ \hline
$10_{162}$&$-3\,0:-2\,0:-2\,0$&$\{\{10\},\{6,10,14,18,16,4,-20,2,8,-12\}\}$&$\{6,-3,2,-1,8,-9,-4,5,10,-7\}$\\ \hline
$10_{163}$&$8^*-3\,0$&$\{\{10\},\{6,-10,14,16,-4,18,2,20,12,8\}\}$&$\{2,7,-6,-3,10,-1,8,-5,4,-9\}$\\ \hline
\end{tabular}

\normalsize

\bigskip

For knots that don't have minimal OGC diagrams, we can search for their non-minimal OGC diagrams. Notice that some knots that can be represented only by a non-minimal meander diagram can have non-minimal OGC diagram with a smaller number of crossings. E.g., knot $9_{49}=-2\,0:-2\,0:-2\,0$ has a meander diagram given by the short Gauss code $\{7,-4,3,-8,9,2,5,-6,-1,-10,11\}$ with $n=11$ crossings, and non-minimal OGC diagram $\{-2,5,-6,1,-10,-3,4,$ $-9,8,-7\}$ with $n=10$ crossings. We have not succeeded to find meander diagrams for only three knots with $n=10$ crossings, $10_{96}=.2.2\,1.2$, $10_{99}=.2.2.2\,0.2\,0$, and $10_{123}=10^*$. Knowing that knot $10_{99}$ has a minimal OGC diagram given by the short Gauss code $\{-4,5,-8,1,-2,9,-10,3,-6,7\}$, that knots $10_{96}$ has non-minimal OGC diagram given by the short Gauss code $\{-6,5,-12,7,2,-1,8,11,-4,3,-10,9\}$ with $n=12$ crossings, and that knot $10_{123}$ has  non-minimal OGC diagram given by the short Gauss code $\{-6,1,-8,11,$ $4,-3,10,-9,-2,5,-12,7\}$ with $n=12$ crossings, we have the following theorem:

\begin{theorem}
Every knot with at most $n=10$ crossings has an OGC diagram.
\end{theorem}

\begin{conjecture}
Every rational knot is OGC knot.
\end{conjecture}

This conjecture is checked for all rational knots with at most $n=15$ crossings.

A diagram is positive if all its crossings are of the same sign. A knot is called {\it positive} if it has a positive diagram. Every alternating positive knot has a minimal positive diagram \cite{nak}. In the paper \cite{stoi} A.~Stoimenov found a positive knot $K11n{183}=9^*.-2:.-2$ without positive minimal diagrams. This knot has only one minimal diagram, and it is not positive. However, its has the positive non-minimal 12-crossing diagram $\{\{12\},\{6,14,10,24,20,4,18,2,22,12,8,16\}\}$.

\begin{theorem}
Every positive OGC diagram with an odd number of crossings is a meander diagram.
\end{theorem}

\noindent {\bf Question:} Is it true that every knot has an OGC diagram?

\subsection{Product of OGC knots}

\begin{definition}
A graph consisting of a circle and $n$ chords joining $2n$
different points on it is called {\it chord diagram} of order $n$.
\end{definition}

Planar diagrams of oriented knots are characterized by their Gauss
diagrams. A Gauss diagram of a classical knot projection is an
oriented circle considered as the preimage of the immersed circle
with chords connecting the preimages of each crossing. To maintain
information about overcrossings and undercrossings, the chords are
oriented toward the undercrossing point and can be also signed in accordance
with crossing signs.

\begin{figure}[th]
\centerline{\psfig{file=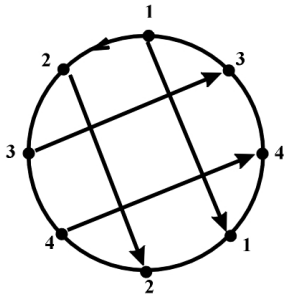,width=1.6in}} \vspace*{8pt}
\caption{Chord diagram of figure-eight knot.
\label{chord}}
\end{figure}

Fig. \ref{chord} shows the chord diagram of a figure-eight knot given by the Gauss code $\{1,-2,3,-4,2,-1,4,-3\}$. Following the orientation of the circle we draw the sequence of numbers $1,-2,3,-4,2,-1,4,-3$ and join $(1,-1)$, $(2,-2)$, $(3,-3)$, and $(4,-4)$ by oriented chords.

Let be given two alternating OGC knots with the same number of crossings $n$ with their ordered Gauss codes $G_1$ and $G_2$. If in $G_1$ and $G_2$ every first half of the code is the sequence of $n$ numbers beginning with $1$ and the absolute value of this sequence is $1,2,\ldots ,n$, $G_1'$ and $G_2'$ are the second halves of the codes, the new Gauss code obtained as concatenation of $G_1'$ and $-G_2'$ we call the {\it product} of ordered Gauss codes $G_1$ and $G_2$.

\begin{theorem}
Product of ordered Gauss codes can be transformed into an ordered Gauss code.
\end{theorem}

\begin{figure}[th]
\centerline{\psfig{file=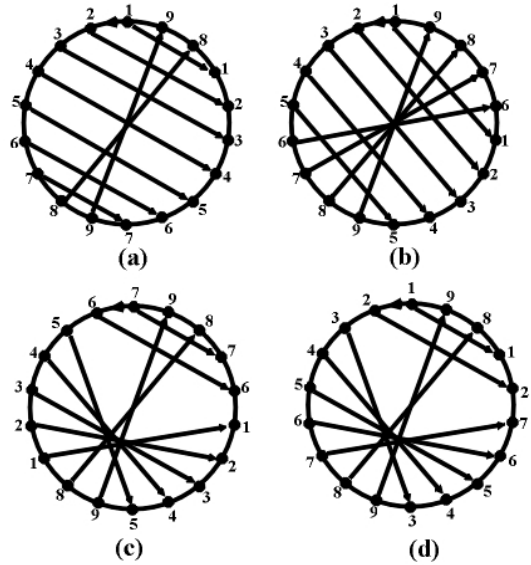,width=2.8in}} \vspace*{8pt}
\caption{Chord diagrams of (a) knot $9_2$; (b) knot $9_4$; (c) their product, knot $9_6$; (d) chord diagram of the ordered Gauss code of the knot $9_6$.
\label{ch}}
\end{figure}

This statement follows immediately from the fact that the first half of the code obtained by the concatenation contains $n$ different numbers.  Such a product of ordered Gauss codes is well defined only for knots with the same odd number of crossings, and for knots with an even number of crossings it results in non-realizable Gauss codes. E.g., from ordered Gauss codes $\{1,-2,3,-4,5,-6,7,-8,$ $9,-7,6,-5,4,-3,2,-1,8,-9\}$ and $\{1,-2,3,-4,5,-6,7,-8,9,-5,4,-3,2,-1,$ $6,-7,8,-9\}$ of knots $9_2=7\,2$ and $9_4=5\,4$ we obtain the code of the product $\{-7,6,-5,4,-3,2,-1,8,-9,5,-4,3,-2,1,-6,7,-8,9\}$, corresponding to the knot $9_6=5\,2\,2$. By using chord diagrams of the knots, the product can be modeled if we take second halves of chord diagrams, make from them an oriented circle and join the corresponding points defining chords (Fig. \ref{ch}). It is clear that the product of a knot with an odd number of crossings given by its ordered Gauss code and its mirror image represents an unknot. Product of every knot with an odd number of crossings $2n+1$ with itself is a torus knot $T(2n+1,2)$, i.e., knot $3_1$, $5_1$, $7_1$, $9_1$, $\ldots $.

For ordered Gauss codes $G_1$ and $G_2$ of knots with the same even number of crossings we can define another kind of product, resulting in a meander link given by the Gauss code $\{\{G_1'\},\{G_2'\}\}$. E.g., the product of Gauss codes $\{1,-2,3,-4,5,-6,7,-8,6,-5,4,-3,2,-1,8,-7\}$ and $\{1,-2,3,-4,5,-6,7,-8,$ $4,-3,2,-1,8,-7,6,-5\}$ of knots $8_1=6\,2$ and $8_3=4\,4$ is the Gauss code $\{\{6,-5,4,-3,2,-1,8,-7\},\{-4,3,-2,1,-8,7,-6,5\}\}$ of the meander link $8_1^2=8$.

In a similar way we can define a product of two 2-component meander links $G_1$ and $G_2$ with the same number of crossings $n$ as a Gauss code $\{\{G_1'\},\{G_2'\}\}$ consisting from the second parts of their codes, i.e., their short Gauss codes. E.g., the product of 12-crossing meander links $7\,5$ given by Gauss code $\{\{-1,2,-3,$ $4,-5,6,-7,8,-9,10,-11,12\},\{1,-6,5,-4,3,-2,7,-8,9,-10,11,-12\}\}$ and \newline $6^*4.2:2\,0.2$ given by Gauss code $\{\{-1,2,-3,4,-5,6,-7,8,-9,10,-11,12\},\{-2,$ $7,-6,3,-12,11,-10,9,-4,5,-8,1\}\}$ is the meander link given by Gauss code $\{\{1,-6,5,-4,3,-2,7,-8,9,-10,11,-12\},\{2,-7,6,-3,12,-11,10,-9,4,-5,$ \newline $8,-1\}\}$ that can be ordered as $\{\{1,-2,3,-4,5,-6,7,-8,9,-10,11,-12\},\{6,-7,$ $2,-5,12,-11,10,-9,4,-3,8,-1\}\}$.

\end{document}